\documentclass[12pt]{article}

\usepackage{amsmath,amsfonts,cd2}
\usepackage{amssymb}

\newtheorem{satz}{Theorem}[section]

\newtheorem{kor}[satz]{Corollary}
\newtheorem{lem}[satz]{Lemma}

\newtheorem{bem}[satz]{Remark}

\newtheorem{prop}[satz]{Proposition}

\newtheorem{pbm}[satz]{Problem}

\newcounter{Roma}
\newcounter{Ara}
\newcounter{let}

\makeindex

\begin{document} 






\newcommand{\nc}{\newcommand}

\nc{\mapco}{\,\colon\, }
\nc{\ab}{^{ab}}


\nc{\hra}{\hookrightarrow}

\nc{\epi}{epimorphism}
\nc{\repi}{regular epimorphism}
\nc{\mono}{monomorphism}
\nc{\iso}{isomorphism}
\nc{\coker}[1]{\mbox{${\rm Coker}(#1)$}}
\nc{\Ker}[1]{\mbox{Ker$(#1)$}}
\nc{\defgl}{\stackrel{def}{=}}

\nc{\V}{\vspace{3mm}}
\nc{\VV}{\vspace{4mm}}
\nc{\lra}{\longrightarrow}
\nc{\lla}{\longleftarrow}
\nc{\mr}[1]{ \stackrel{#1}{\lra} }
\nc{\ml}[1]{ \stackrel{#1}{\lla} }
\nc{\hmr}[1]{\hspace{2mm}\stackrel{#1}{\lra}\hspace{2mm}}
\nc{\hml}[1]{\hspace{2mm} \stackrel{#1}{\lla}\hspace{2mm}}
\nc{\N}{\noindent}
\nc{\st}{^{\prime}}
\nc{\ot}{\otimes}
\nc{\hcong}{ \hspace{2mm}\cong\hspace{2mm}  }
\nc{\hfbox}{\hfill$\Box$}
\nc{\REF}[1]{(\ref{#1})}

\def\Z{\ifmmode{Z\hskip -4.8pt Z} \else{\hbox{$Z\hskip -4.8pt Z$}}\fi}

\def\Q{\ifmmode{Q\hskip-5.0pt\vrule height6.0pt depth
0pt\hskip6pt}\else{\hbox{$Q\hskip-5.0pt\vrule height6.0pt depth
0pt\hskip6pt$}}\fi}

\newcommand{\NN}{\mbox{$I\!\!N$}}

\nc{\Ph}{\phantom{}}

\nc{\BE}{\begin{equation}}
\nc{\EE}{\end{equation}}

\nc{\dst}{\displaystyle}
\nc{\sst}{\scriptscriptstyle}
\nc{\ssst}{\scriptscriptstyle}
\nc{\proof}{\N{\bf Proof\,:}\quad}
\nc{\proofof}[1]{\N{\bf Proof of {#1}\,:}\quad}
\nc{\proofofthm}[1]{\N{\bf Proof of Theorem \ref{#1}\,:}\quad}

\nc{\htt}[1]{^{\otimes #1}}

\newcommand{\ssur}[2]{\mbox{$#1 \!\to\!\!\!\!\!\to\! #2$}}

\def\mapsel#1{\mbox{$\rule[-5mm]{0mm}{12mm} \searrow\rlap
{$\vcenter{\makebox[0mm][r]{$\scriptstyle#1$\hspace{12.2mm}}}$}$}}

\def\surltop#1{\makebox[1cm]{\mbox{$\stackrel{#1}{\makebox[0mm]{$\lla$}\hspace{0.7mm}
\makebox[0mm]{$\lla$}}$}}}

\nc{\Sur}[1]{\mbox{$\:\stackrel{#1}{\lra\!\!\!\!\!\to\,}\:$}}

\def\INJ{\mbox{\mathsurround=0pt
\makebox[0mm][r]{\parbox{0mm}{\rule[-0.65mm]{0mm}{0.2mm}$\scriptscriptstyle>$}}
\makebox[0.7cm][l]{\parbox{0.7cm} {$\lra$}}}}

\def\Inj#1{\mbox{$\:\stackrel{#1}{\INJ}\:$}}

\nc{\Injup}[1]{\mapup{#1}}
\nc{\injup}[1]{\mapup{#1}}
\nc{\injdown}[1]{\mapdown{#1}}


\def\mapup#1{\mbox{$\rule[-1mm]{0cm}{0.7cm}  
\makebox[0mm][r]{\raisebox{0.2mm}{$\scriptstyle\phantom{\cong}$}\hspace{0.6mm}}
\bigg\uparrow\rlap{$\vcenter{\hbox{$\scriptstyle#1$}}$}$}}

\def\isoup#1{\mbox{$\rule[-1mm]{0cm}{0.7cm}  
\makebox[0mm][r]{\raisebox{0.2mm}{$\scriptstyle\cong$}\hspace{0.6mm}}
\bigg\uparrow\rlap{$\vcenter{\hbox{$\scriptstyle#1$}}$}$}}

\def\mapdown#1{\mbox{$\rule[-1mm]{0cm}{0.7cm}  
\makebox[0mm][r]{\raisebox{0.2mm}{$\scriptstyle $}\hspace{0.6mm}}
\bigg\downarrow\rlap{$\vcenter{\hbox{$\scriptstyle#1$}}$}$}}

\def\isodown#1{\mbox{$\rule[-1mm]{0cm}{0.7cm}  
\makebox[0mm][r]{\raisebox{0.2mm}{$\scriptstyle\cong$}\hspace{0.6mm}}
\bigg\downarrow\rlap{$\vcenter{\hbox{$\scriptstyle#1$}}$}$}}

\def\isor#1{\mbox{$\smash{\mathop{\longrightarrow}\limits^{\cong}_{#1}}$}}

\def\isol#1{\mbox{$\smash{\mathop{\longleftarrow}\limits^{\cong}_{#1}}$}}

\def\surdown#1{\makebox[0mm]{$\mapdown{\raisebox{0.4mm}{$\scriptstyle#1$}}$}
\makebox[0mm]{\raisebox{-2.15mm}{$\downarrow$}}}

\def\surup#1{\makebox[0mm]{$\mapup{\raisebox{0.4mm}{$\scriptstyle#1$}}$}
\makebox[0mm]{\raisebox{0.7mm}{$\mapup{}$}}}

\newcommand{\surr}[2]{@\cdhgeneric>->->\twoheadrightarrow>#1>#2>}

\newcommand{\surl}[2]{@\cdhgeneric>\twoheadleftarrow>->->#1>#2>}

\newcommand{\brokenr}[2]{
@\cdhgeneric>\raise2.8pt\hbox to3.5pt{\hrulefill}\mkern9mu>
\raise2.8pt\hbox to3.5pt{\hrulefill}\hbox to5pt{}>
\mkern-7.5mu\dashrightarrow>#1 >
#2>}

\newcommand{\brokenup}[1]{
@\cdvgeneric>\hat\cdot
>\raisebox{3pt}{$\vdots$}>
\vbox{\kern3.5pt\hbox{$\cdot$}\kern-3.5pt}
> >\hspace{1mm}#1>}

\newcommand{\functlr}[2]{
\raisebox{0.4pt}{$\hss\begin{CD}
@>\vbox{\hbox to 0pt{$\hss\begin{CD}@<#1<<\end{CD}\hss$}\vskip-2pt}
>#2 >
\end{CD}\hss$}
}

\newcommand{\maprr}[2]{
\raisebox{-0.9pt}{$\hss\begin{CD}
@>\vbox{\hbox to 0pt{$\hss\begin{CD}@>#1>>\end{CD}\hss$}\vskip-3pt}
>#2 >
\end{CD}\hss$}
}

\newcommand{\mapdd}[2]{
@\cdvstandard>\downarrow\hspace{1.5pt}\kern-4pt\hspace{1.5pt}\Big\downarrow>#1>#2>}

\newcommand{\mapud}[2]{
@\cdvstandard>\uparrow\kern-3pt\Big\downarrow>#1>#2>}

\newcommand{\sepi}[3]{\,\mbox{$#1\,: \ssur{#2}{#3}$}\,}

\nc{\auf}{\twoheadrightarrow}
 
\nc{\ruled}{\rule[-4mm]{0mm}{0mm}}

\nc{\sm}{\:{ \wedge}\:}


\nc{\map}[3]{\mbox{$#1 \mapco #2 \to #3$}}

\nc{\rond}{{\,\sst \circ\,}}

\nc{\ruleu}{\rule{0mm}{7mm}}

\nc{\T}[1]{\tilde{#1}}

\nc{\Imm}[1]{\mbox{${\rm Im}(#1)$}}

\nc{\tw}{\end{document}}

\nc{\QJG}{\frac{\dst I(G) J}{\dst I^2(G) J}}
\nc{\QJH}{\frac{\dst I(H) J}{\dst I^2(H) J}}
\nc{\otz}{\ot}
\nc{\UL}[2]{{\rm U}_{#1}{\rm L}(#2)}

\nc{\IRN}[1]{I_{R,\cal G}^{#1}(G)}

\nc{\ULG}[1]{{\rm U}_{#1}{\rm L}^{\cal G}(G)}
\nc{\UGH}[1]{{\rm U}_{#1}^{\cal G}(G,H)}

\nc{\LG}[1]{{\rm L}^{\cal G}_{#1}(G)}

\nc{\calG}{{\cal G}}
\nc{\AB}{^{AB}}

\vspace*{0.7cm}

\begin{center} \LARGE \bf On Fox quotients of arbitrary group algebras \normalsize  
 \vspace{0.8cm}

\author{MH}

\renewcommand{\thefootnote}{ }
\normalsize \bf  Manfred Hartl\normalsize \footnote[1]{Research partially supported by an individual
fellowship of the Human Capital and
Mobility }
\footnote[1]{Programme of the European Union during the years 1993-1995.}
\footnote[1]{2000 {\em Mathematical subjects classification}\/  20C07, 20J05, 17B35.}
\vspace{5mm}

\N \rm Universit\'e de Lille - Nord de France,\\
LAMAV  and FR CNRS 2956,\\ 
ISTV2, UVHC\\ 
Le Mont Houy,  59313 Valenciennes Cedex 9, France.\\ 
Email: Manfred.Hartl@univ-valenciennes.fr\\ 

\end{center}
\vspace{8mm}

\begin{abstract}

\N Certain subquotients of group algebras are determined as a basis for subsequent computations of relative Fox and dimension subgroups. More precisely,
for a group $G$ and  N-series ${\cal G}$ of $G$ let $I^n_{R,\cal G}(G)$, $n\ge 0$, denote the filtration of the group
algebra $R(G)$ induced by ${\cal G}$, and $I_R(G)$ its augmentation ideal.
For subgroups $H$ of $G$, left ideals $J$ of $R(H)$ and right $H$-submodules $M$ of $I_{Z\!\!\!Z}(G)$
the quotients $I_R(G)J/MJ$ are studied by homological methods, notably for
 $M= I_R(G)I_R(H)$, $I_R(H)I_R(G) + I([H,G])R(G)$ and
$R(G)I_R(N) +I_{R,\cal G}^n(G)$ for a normal subgroup $N$ in $G$; in the latter case the module  $I_R(G)J/MJ$ is completely determined for $n=2$.
The groups  $I_{Z\!\!\!Z,\cal G}^{n-1}(G)I_{Z\!\!\!Z}(H)/I_{Z\!\!\!Z,\cal G}^n(G)I_{Z\!\!\!Z}(H)$
  are studied and explicitly computed for $n\le 3$  in terms of enveloping rings of certain graded Lie rings
and of torsion products of abelian groups.
\end{abstract}

\N{\bf Keywords\,:} group algebra, augmentation quotient, Fox subgroup, N-series, enveloping algebra.\vspace{5mm}

\begin{center} {\bf Introduction}\V \end{center}

Let $G$ be a group, $R(G)$ be the group algebra of $G$ with coefficients in a commutative ring $R$,  and let $I^n_{R,\cal G}(G)$,
$n\ge 0$, denote the filtration of $R(G)$ induced by a given N-series ${\cal G}$ of $G$, see \cite{Pa} or section
2 below. In particular,   $I^0_{R,\cal G}(G) = R(G)$, $I^1_{R,\cal G}(G)$ is the augmentation
ideal $I_R(G)$ of $R(G)$, and if $\calG$ is the lower central series $\gamma$ of $G$, then $I^n_{R,\gamma}(G)$ is the $n$-th power
$I_R^n(G)$ of $I_R(G)$. As usual, we skip the sub- or superscript $R,\calG$ when $R=\Z$ or $\calG=\gamma$, resp. 

 Now let $H$ be a subgroup of $G$ and $J$ be a left ideal of $R(H)$. Then it
is a \textit{classical problem} to study the \textit{$n$-th generalized Fox quotient} 
$ I_{R,\cal G}^{n-1}(G) J/ I_{R,\cal G}^n(G) J$  for $n \ge 1$ (usually one considers $R=\Z$, $\calG =\gamma$ and $J= I(H)$). The
first Fox quotient was determined a long time ago              
\cite{Whitcomb},  and the related quotient  $ I(G) I(H)^{n-1}/ I(G) I(H)^n$ was 
computed by Karan and Vermani \cite{Ka-Ve}. But only partial results are known on the second
Fox quotient (and some higher ones),  for  $R=\Z$  and always under some splitting condition, e.g.\ supposing $H$  to be a
semi-direct factor of $G$, see the work of Khambadkone and of Karan and Vermani in
\cite{Kh12}, \cite{Kh13}, \cite{Kh14}, \cite{KV29},
\cite{KV30}.  
As a striking fact,  in all known cases a reduction of the  
quotient in question was obtained to \textit{sums and tensor products} among certain subquotients of the group rings
of $G$ and $H$  
(or of a semidirect complement of $H$ in $G$):  examples of this phenomenon are given by Proposition \ref{Fox2split} and by the isomorphisms \REF{introKa-Ve} and \REF{introKVgen} below.

 
In this paper we introduce a
\textit{homological approach} to the study of the quotient $I_R(G)J/MJ$
for arbitrary $G$, $H$, $R$, $J$ and right $H$-submodules $M$ of $I(G)$, in the same spirit: 
in certain cases we still find a sum-tensor-decomposition, such as the isomorphisms
\BE\label{introKa-Ve}
\mbox{$
 \frac{ I_R(G)J}{  I_R(G)I_R(H)J} \hspace{2mm} \cong\hspace{2mm} \frac{  I(H)J}{ 
I^2(H)J} \hspace{2mm}\oplus\hspace{2mm} \big (\frac{  I(G) }{ 
\Z(G)I(H) }\big ) \ot \big ( \frac{  J}{  I(H)J} \big )
$}
\EE
in Theorem \ref{Ka-Ve} and
\BE\label{introKVgen}\mbox{$
\frac{ I_R(G) J}{ I_R(H)I_R(G) J +  I_R([H,G]) R(G)J}   \hspace{2mm}\cong\hspace{2mm}  
\big(\frac{ H}{ [H,G]}  \hspace{0.5mm}\ot\hspace{0.5mm} \frac{ 
J}{  I(H)J}\big)\hspace{2mm}\oplus\hspace{2mm}  \big(I(G/H) \hspace{0.5mm}\ot\hspace{0.5mm} \frac{ 
J}{ I(H)J}\big)\ruled$}
\EE
in Theorem \ref{KVgen}, the latter for normal subgroups $H$. But for general submodules $M$, the group
 $I_R(G)J/M_RJ$ can only be embedded
into a natural exact sequence
\BE\label{introgenFoxsequ}\mbox{$
{\rm Tor}_1^H(\frac{I(G)}{I(H)+M}\,,  J)  
\hspace{1mm}\mr{\sigma} \hspace{1mm} 
\frac{\rule{0mm}{4mm}  I(H) J  }{   (I(H)   \cap M) J }     \hspace{1mm}\mr{j} \hspace{1mm} 
  \frac{  I_R(G) J}{  M_R J   \rule{0mm}{2.7mm}\ruled} 
\hspace{1mm}\mr{} \hspace{1mm}
  \frac{ I(G)}{I(H)+M} 
\hspace{1mm} \ot_H \hspace{1mm}   J  \hspace{1mm} \mr{}\hspace{1mm} 0$}
\EE
\N 
thus  describing the kernel and cokernel of the canonical map  $j$  in terms of a \textit{tensor and a torsion
product over
$H$}, see Proposition \ref{extforM}. Moreover, the class of  the resulting  extension $0\to \coker{\sigma} \to I_R(G) J/ M_R J \to (I(G)/I(H)+M) 
\hspace{1mm} \ot_H \hspace{1mm}   J \to 0$ is also described, in terms of another, more accessible extension.


If $R=\Z$ and $J =I(H)$  sequence \REF{introgenFoxsequ} admits a long exact extension, thus describing \Ker{j} in terms of  \textit{homology groups} of
$H$ (Theorem \ref{H2sequ}). This also provides a   reduction of  \textit{Fox subgroups of $G$} to \textit{induced subgroups of $H$}, see Corollary \ref{Foxsgp}; for example if $H$ is \textit{free} then  
\BE\label{introFoxsgr} G\cap (1+MI(H))= H\cap (1+(M\cap I(H))I(H))\,.\EE

A case of particular interest is  the quotient
$I_R(G)J/ (R(G)I_R(N)J+I_{R,\calG}^n(G)J) $ for a normal subgroup $N$ of $G$; when $J=I_R(H)$ we call it the \textit{Fox polynomial group relative to  $N$ and $\calG$}
since it generalizes the  \textit{polynomial
group}
$P_{n,R}(G) = I_R(G)/ I^{n+1}_R(G)$ introduced by Passi in \cite{Pa68}, and its relative version
$P_{n,R}(G,N) = I_R(G)/ (I_R(N)I_R(G)+I^{n+1}_R(G))$  used in various contexts, cf.\ section 2.

Just as polynomial groups serve to analyze 
augmentation quotients and, dually, dimension subgroups, we  here show how the study of Fox polynomial groups leads to
a computation of certain Fox quotients  and Fox subgroups.  The descriptions obtained formally generalize known results in the case where $H=G$, see Theorems \ref{Fox2IH} and \ref{Fox3push}; they reveal  a striking duality between (generalized) Fox polynomial groups and (generalized) Fox subgroups, which comes from the symmetry of push-outs. Fox subgroups are further studied
in subsequent work, e.g.\ in \cite{D3F2}; in this paper  we focus on
Fox polynomial groups and Fox quotients which we compute in the following
cases: 
If the group $H$ is \textit{abelian}, we determine the \textit{relative Fox polynomial group} 
$I(G)I(H)/(\Z(G)I(N)I(H)+I^n(G)I(H))$ for all
$n$ and normal subgroups $N$ of $G$, in Corollary \ref{Hfreeabcyc}.
The \textit{first generalized Fox quotient} is given by
\[\mbox{$
\frac{ R(G)J}{I_R(G)J} \hspace{2mm} \cong \hspace{2mm} \frac{ J}{I(H)J}
$}\]
see Proposition \ref{Fox1} which generalizes a well-known result of Whitcomb for $J=I(H)$. 
The \textit{second generalized Fox quotient} is determined in Theorem \ref{Fox2}, by means of a natural exact sequence, for any subgroup $K$ of $G$:\V

\N\makebox[14.7cm]{ \makebox[0mm]{
\begin{minipage}{16cm} \small
\[\begin{matrix}
\stackrel{}{{\rm Tor}_1^{\Z}(\frac{G}{HKG_{(2)}}\,,\frac{J}{I(H)J})} &   & & &  & &\cr
\downarrow & & & &  & &\cr
\frac{  I_R(H)J}{  I_R(H \cap KG_{(2)}) J+ I^2_{R}(H) J} & \mr{j }  & 
\frac{   I_R(G)J}{  I_R(K) J + I_{R,{\cal G}}^2(G) J}
  & \mr{}  &  \big( \frac{  G}{ HKG_{(2)}}\big)
\ot  \big(\frac{J}{I(H)J\rule[-6mm]{0mm}{0mm}} \big) & \mr{}  & 0
 \end{matrix}\]
\end{minipage}
}}
A case where the bottom sequence is split short exact is given in Corollary \ref{Fox2split}; this in particular covers the special cases treated in the literature. The above sequence also allows to explicitly compute the intersection $I_R(H)J \cap
(I_R(K)J+I^2_R(G)J)$, see Corollary \ref{intersec}.

As to the higher Fox quotients $Q_{n,R}^{\cal G}(G,H)=I^{n-1}_{R,\calG}(G)I_R(H)/I^n_{R,\calG}(G)I_R(H)$,
Proposition \ref{Hfreequot} provides a natural isomorphism
\BE\label{introQnHfree}\mbox{$
Q_{n,R}^{\cal G}(G,H)
\hspace{2mm} \cong \hspace{2mm} Q_{n,R}^{\cal G}(G,1\!\!1) \hspace{1mm}\ot\hspace{1mm} H^{ab}$}\EE
%
if $H$ is a \textit{free group}; if only its $n-1$-step nilpotent quotient is free nilpotent then the right-hand term in \REF{introQnHfree} must be replaced by a natural quotient, see Theorem \ref{Foxfreequot}. The resulting description of $Q_{n}(G,H)
$ in this case strongly resembles
 the description of $Q_{2}(G,H)$
 for  arbitrary subgroups $H$ given in Theorem \ref{Fox2IH}. 


In order to study the higher Fox quotients   $Q_n^{\cal G}(G,H)$ for arbitrary subgroups $H$, we combine our homological methods with a suitable \textit{generalization of Quillen's approximation of
augmentation quotients} \cite{Quillen}. Recall from \cite{Pa} that for any N-series $\calG$ of $G$ there is a surjective morphism of graded rings
\[ \mbox{$\theta^{\calG}\,\colon\,\ULG{} \hspace{2mm} \Sur{} \hspace{2mm}
 \mbox{Gr}^{\calG}(\Z(G)) =
\sum_{n\ge 0} I^n_{\cal G}(G)/I^{n+1}_{\cal G}(G)$}\]
 where $\ULG{}$ is the enveloping algebra of the graded Lie ring $\LG{}$ formed
by the successive
 quotients of $\calG$. As a key fact, Quillen  showes in  \cite{Quillen} that rationally $\theta^{\gamma}$ is an isomorphism, i.e.
\Ker{\theta^{\gamma}} is torsion; this is in fact still true for any N-series as follows from  fundamental work of Hartley
\cite{Hartley}, see also
\cite{Pa}. We here extend  this construction to Fox quotients by introducing  a surjective graded bimodule map
\[\mbox{$ \theta^{\calG H} \,\colon\, \ULG{}
\ot_{{\rm UL}(H)}  {\rm
\bar{U}L}(H)
\hspace{2mm} \Sur{} \hspace{2mm}
\sum_{n\ge 1} I^{n-1}_{\cal G}(G)I(H)/I^{n}_{\cal G}(G)I(H)$}\]
 which coincides with
$\theta^{\calG}$ if $H=G$ and $\calG=\gamma$. It turns out that just like
$\theta^{\calG}$ \cite{Q3},
\textit{$\theta^{\calG H}$ is an isomorphism in degrees 1 and 2} (Proposition \ref{xi1et2}) but not in higher degrees, in general. We exhibit canonical elements in its
kernel (Proposition \ref{kerneltheta}); this leads to the problem of whether they generate
\Ker{\theta^{\calG H}} modulo torsion, in extension of Quillen's theorem. In degree 3 this is shown to hold, based on the natural exact sequence
 \[ \ruleu {\rm Tor}_1^{\ifmmode{Z\hskip -4.8pt Z} \else{\hbox{$Z\hskip -4.8pt Z$}}\fi}
 (G^{AB},   H\ab) \:\oplus\: \Big(\Ker{l_2^{\calG H}}\cap \Ker{c^H_2} \Big)
\:\mr{(\delta_1,\delta_2)}\: \UGH{3} \:\mr{\theta_3^{\calG H}} \:Q_3^{\calG}(G,H) \:\to\: 0 
\ruled \] 
established in Theorem \ref{Fox3}. This
amounts to a \textit{complete description} of  the structure of $I^2 _{\cal G}(G)I(H)/I^3 _{\cal G}(G)I(H)$, both   in a functorial way and in terms of a cyclic decomposition of $H\ab$ if $H$ is finitely
generated. The result obtained resembles our description of
\Ker{\theta_3^{\calG}}  in \cite{Q3}; the situation here, however, is   more intricate as it involves a
\textit{secondary operator} (namely $\delta_2$ which is an additive relation with indetermancy ${\rm Im}(\delta_1)$), a phenomenon which does not occur in the computation of \Ker{\theta_3^{\calG}}. 

The resulting, in fact less complicated  description of the group $Q_3^{\calG}(G,G )$ in Theorem \ref{Q3(G,G)} and Corollary \ref{Q3(G,G)formel} also
 amounts to a functorial computation  of the group 
\[
\frac{I^3(H) \oplus I([H,K])I(H) }{ I^4(H) + I(H)I([H,K])I(H) + I([H,K,H])I(H) + I([H,K,K])I(H)}\]
taking $H=G$ and a specific N-series ${\calG}$; this group
 appears in Khambadkone's analysis  of Fox quotients of semidirect products \cite{Kh12}, but seems  not to have been computed yet, not  even  in   special cases.

Finally, we point out that the description of \Ker{\theta_3^{\calG G}} above also allows to determine the \textit{ fourth relative dimension subgroup} $G\cap (1+I(N)I(G)+I^4(G))$ for arbitrary normal subgroups $N$ of $G$ (work in progress).

\vspace{4mm}

\N{\bf Conventions.}\quad  In this paper, the terms of the lower central series of $G$ are denoted by $G_i=[G_{i-1},G]$, with
$G_1=G$ and $[a,b] = aba^{-1}b^{-1}$ for $a,b\in G$. We write $G\ab = G/G_2$. 
Maps denoted by $i$ or $j$ (possibly indexed) are
always induced by the inclusion $H\hra G$ or $I(H) \hra I(G)$, and 
 $q$ 
(possibly indexed) always
denotes a canonical quotient map. Moreover, maps denoted by $\mu_G,\mu_G\st$ or $\mu_H,\mu_H\st$ are always
induced by multiplication in $\Z(G)$ or $\Z(H)$, respectively, while 
$\epsilon\,\colon\,R(G) \to R$ denotes  the augmentation map.
 Maps denoted  by $\bar{f}$ defined on some quotient
$A/B$ are induced by a homomorphism $f$ on $A$. Tensor products over the ring $\Z(H)$ are denoted by $\ot_H$ and over $\Z$ by $\ot$.

%


\tableofcontents

\section{Exact sequences for Fox quotients}

In this section we provide the key tools (the ``pushout lemma'' \ref{pushlem} and the generic exact sequences \ref{pushsequ} and
\ref{extforM}) of our  study of the Fox problem in various generalized forms. In particular it allows to deal with
Fox quotients and Fox subgroups simultaneously and in a strictly dual way, cf.\ Theorems \ref{Fox2IH} and \ref{Fox3push}. For
general properties of pushouts see \cite{Rotman}.  

Recall the notation from the introduction. We point out that we throughout use the plain facts that $I_R^n(G)J=I^n(G)J$ for
$n\ge 0$, and that $I_R(G) \cong I(G) \ot R$, whence for any subgroup $M$ of $I(G)$,
  \[ \frac{\dst I(G)}{\dst M} \ot J  \hspace{2mm}\cong \hspace{2mm}    \frac{\dst I(G)\ot J}{\dst \Imm{M\ot J}} 
  \hspace{2mm}\cong \hspace{2mm}   \frac{\dst I_R(G)\ot_R J}{\dst \Imm{M_R\ot_R J}}   
  \hspace{2mm}\cong \hspace{2mm}   \frac{\dst I_R(G)}{\dst M_R} \ot_R J \]
where $M_R$ is the $R$-submodule of $I_R(G)$ generated by $M$. This also implies that 
$\frac{I(G)}{ M} \ot_H J  \hspace{2mm}\cong \hspace{2mm}  \frac{I_R(G)}{ M_R} \ot_{R(H)} J$.


 We start with 
the following elementary

\begin{lem}\label{lem1}\quad The canonical map
  \[ \sepi{\mu_J}{\Z(G) \ot_H J\:}{\:\Z(G)  J = R(G)  J}  \]
is an isomorphism.
\end{lem}

\N{\bf Proof\,:}\quad Using the injections  $\eta\,\colon\,J \hra R(H)\stackrel{\cong}{\to} \Z(H)\ot R$ and $\eta\st\,\colon\,R(G)J
\hra R(G)$ we have the following commutative square
 \[ \begin{matrix}
\ruleu\Z(G) \ot_H J  &  \Sur{\mu_J} & \Z(G)  J \cr
\mapdown{id \ot \eta}  &  & \injdown{\eta\st} \cr
\ruled\Z(G) \ot_H \Z(H) \ot R & \cong  & R(G)   \ruleu
\end{matrix}\]
 But $id \ot \eta$ is injective since $\Z(G)$ is a free $H$-module, so also $\mu_J$
is injective. $\Box$
\vspace{3mm}

In order to describe our key lemma, let $U \subset M \subset I(G)$ be right $H$-submodules and $V \subset N
\subset J$ be left $R(H)$-submodules such that $U \subset I(H)$. Consider the following
commutative square.
  \BE\label{dia1}
\begin{CD} \begin{matrix}
  \frac{\ruleu\dst I(H) \ot_H J}{\ruled\rule{0mm}{3.5mm}\dst \Imm{U \ot_H V}}   &  \mr{\overline{i\ot id}}  &  
  \frac{\dst I(G) \ot_H J}{\dst \rule{0mm}{3.5mm} \Imm{M \ot_H N}}   \cr
\surdown{\bar{\mu}_H}  &  &  \surdown{\bar{\mu}_G}  \cr
  \frac{\rule{0mm}{6mm}\dst I(H)   J}{\dst\rule{0mm}{3.5mm} U   V}  &  \mr{j}  &    
\frac{\dst I(G)   J}{\ruled\rule{0mm}{3.5mm}\dst M N}
\end{matrix}   \end{CD}\EE

\begin{lem}\label{pushlem}\quad Diagram \REF{dia1} is a pushout square of $R$-modules,
or equivalently, $\:\Ker{\bar{\mu}_G} = (\overline{i\ot id})\Ker{\bar{\mu}_H}$\,. This implies that also $\Ker{j} =
\bar{\mu}_H\Ker{\overline{i\ot id}}$.
Moreover, if $J = I_R(H)$ then $\Ker{\bar{\mu}_H} = \pi_H {\rm Im}(\lambda\,\colon\, H_2(H,R) \to I(H) \ot_H I_R(H))$ where $
\pi_H$ is the canonical projection and $\lambda$ is defined in the proof below.
\end{lem}\V

\N{\bf Proof of Lemma \ref{pushlem} and Proposition \ref{Fox1} below\,:}\quad Consider the following
commutative diagram whose middle and bottom row (omitting the right hand bottom corner) are part of the long exact sequences
obtained from tensoring with $J$ the short exact sequences $I(K) \stackrel{\nu_K}{\hra} \Z(K) \Sur{\epsilon_K} \Z$, for $K=G$ resp.\
$H$.
  \[ \begin{CD}\begin{matrix} 
 & & I(G) J &  \hra  &  \Z(G) J  &  \Sur{}  &  \frac{\rule{0mm}{7mm}\dst \Z(G)J}{\dst I(G)J }
\cr
 & & \surup{\mu_G}  &  &  \isoup{\mu_J}  &  & \mapup{\bar{\mu}_J}  \cr
{\rm Tor}_1^H(\Z,J)  &  \mr{\tau_G}  & I(G) \ot_H J  & \mr{\nu_G \ot id}  & \Z(G) \ot_H J 
&  \Sur{\epsilon_G \ot id}  &  \Z  \ot_H  J \cr
\|  & & \mapup{i\ot id} &  & \mapup{i\ot id} &  & \cr
{\rm Tor}_1^H(\Z,J)  &  \mr{\tau_H}  & I(H) \ot_H J  & \mr{\nu_H \ot id}  & \Z(H) \ot_H J
&  \stackrel{\mu_H\st}{\cong} &  J\hspace{9mm} \ruled
\end{matrix}\end{CD}\]
Now $\mu_J$
is isomorphic by Lemma \ref{lem1} and  $\mu_G$ is surjective, so the induced map
$\bar{\mu}_J$ is an isomorphism, too. But $\Z \ot_H J \cong J/ I(H)J $, which proves
Proposition \ref{Fox1}.
Furthermore,
  \begin{eqnarray}\label{East} \Ker{\mu_G} &=& \Ker{\nu \ot id}  = \Imm{\tau_G} = (i \ot id)\Imm{\tau_H} = (i
\ot id) \Ker{\mu_H\st(\nu \ot id)} \nonumber\\
 &=& (i \ot id) \Ker{\mu_H\,\colon I(H) \ot_H J \to I(H)J}\:. 
\end{eqnarray}
Noting $\pi_G\,\colon\,I(G)\ot_H J\auf (I(G)\ot_H J)/\Imm{M\ot_H N}$  the canonical projection we obtain identities 
\begin{eqnarray*}
\Ker{\bar{\mu}_G}  & = & \pi_G(\Ker{\mu_G}) \\
&= & \pi_G(i \ot id)(\Ker{\mu_H})  \\
  &  =  &  (\overline{i \ot id}) \pi_H(\Ker{\mu_H}) \\
  & = & (\overline{i \ot id}) \Ker{\bar{\mu}_H}\:.
\end{eqnarray*}
The identity $\Ker{j} =
\bar{\mu}_H\Ker{\overline{i\ot id}}$ now follows by an easy diagram chase using the surjectivity of $\bar{\mu}_H$, or
more abstractly, by symmetry of  pushouts.
Now consider the case $J=I_R(H)$. First note that for $i\ge 1$, ${\rm Tor}_i^H(\Z,R(H))=0$ since any  resolution $\underline{P}$ of $\Z$ by projective  $\Z(H)$-modules is $\Z$-split exact, whence $\underline{P}\otimes_HR(H) \cong \underline{P}\otimes_H\Z(H)\otimes R\cong \underline{P}\otimes R$ is $R$-split exact.
Then dimension shifting along the short exact sequence $I_R(H) \hra R(H) \Sur{\epsilon} R$ of $H$-modules
provides a connecting isomorphism $\tau\colon H_2(H,R) = {\rm Tor}_2^H(\Z,R) 
\smash{\mathop{\longrightarrow}\limits^{\cong}_{}}
 {\rm Tor}_1^H(\Z,I_R(H))$, so putting 
$\lambda = \tau_H \tau$ also the last part of the assertion is proved. \hfbox

We first consider an easy but useful case which includes the classical case where $G$ is free.

\begin{kor}\label{Hfree}\quad Suppose that $H$ is a free group with basis $X$. Then for any right $H$-submodule $M$ of $I(G)$ there
is an isomorphism  
 \[  I_R(G)I_R(H)/M_RI_R(H) \:\cong\: (I_R(G)/M_R)\ot H\ab  \]
such that for $x\in X$ and $y\in I_R(G)$ the element $(y+M_R) \ot (xH_2)$ corresponds to $y(x-1) + M_RI_R(H)$.
\end{kor}

\proof Let $U=0$ and $J=V=N=I_R(H)$ in \ref{pushlem}. As $H$ is free, $H_2(H,R)=0$ so $\bar{\mu}_G$ is an isomorphism. On the other
hand, $I_R(H)$ is a free $R(H)$-module with
basis
$x-1$,
$x\in X$, see \cite{Hi-St}. This provides a (non natural) isomorphism of left $H$-modules $I_R(H) \:\cong\: \Z(H) \ot H^{ab}\ot R$.
Thus we have isomorphisms 
$I(G)I_R(H)/M I_R(H)   \cong 
I(G)\ot_H I_R(H)/\Imm{M\ot_H I_R(H)} $ $ \:\cong\:(I(G)/M)\ot_H I_R(H)$ $\cong\:
(I(G)/M)\ot H\ab\ot R \:\cong\: (I_R(G)/M_R)\ot H\ab$.\hfbox\V


Using right-exactness of the tensor product Lemma \ref{pushlem} implies the following generic fact which is the basis of all that
follows.\V

\begin{kor}\label{pushsequ} Under the hypothesis of Lemma \ref{pushlem} there is a commutative diagram with exact rows where $M\st
= M+I(H)/I(H)$ and $K=\Ker{\overline{i\ot id}}$.
  \BE\label{pushsequdia}
\begin{CD} \begin{matrix}
K &  \stackrel{\iota}{\hra} &  \frac{\ruleu\dst I(H) \ot_H J}{\ruled\rule{0mm}{3.5mm}\dst \Imm{U \ot_H V}}   &  \mr{\overline{i\ot
id}}  &  
  \frac{\dst I(G) \ot_H J}{\dst \rule{0mm}{3.5mm} \Imm{M \ot_H N}} &
\mr{\overline{q\ot id}}  &  
   \frac{\dst (I(G)/I(H)) \ot_H J}{\dst \rule{0mm}{3.5mm} \Imm{M\st \ot_H N}}  &
\to & 0 \cr
\| & & \surdown{\bar{\mu}_H}  &  &  \surdown{\bar{\mu}_G}  &  &  \| \cr
K &  \mr{\bar{\mu}_H \iota} &   \frac{\rule{0mm}{6mm}\dst I(H)   J}{\dst\rule{0mm}{3.5mm} U   V}  &  \mr{j}  &    
\frac{\dst I(G)   J}{ \ruled\rule{0mm}{3.5mm} \dst M N}  & 
\mr{\overline{\overline{q\ot id}}}  &  
   \frac{\dst (I(G)/I(H)) \ot_H J}{\dst \rule{0mm}{3.5mm} \Imm{M\st  \ot_H N}}  &
\to & 0  
\end{matrix}\end{CD}\EE
Consequently, the bottom sequence is induced by the top one, i.e., 
the abelian group extensions obtained by dividing out the image of $K$ from the left hand terms
are such that the top one induces the bottom one via the map induced by $\bar{\mu}_H$.\hfill$\Box$
\end{kor}

Note that the unhandy terms of the top row of diagram \REF{pushsequdia} reduce to tensor products when either $U=I(H)$ and $M=I(G)$
 or
$V=N=J$. We are mainly concerned with the second case in the rest of the paper, but the first one is also of interest: it allows to
recover a result of Karan and Vermani  
\cite{Ka-Ve} from our viewpoint, in a slightly more general form (arbitrary $R$ and $J$ instead of $\Z$ and $I^{n-1}(H)$).

\begin{satz}\label{Ka-Ve}\quad There is a (non-canonical) isomorphism of $R$-modules\V
  
  \Ph\hfill$  \frac{\dst I_R(G)J}{\dst I_R(G)I_R(H)J} \hspace{2mm} \cong\hspace{2mm} \frac{\dst I(H)J}{\dst
I^2(H)J} \hspace{2mm}\oplus\hspace{2mm} \bigg(\frac{\dst I(G) }{\dst
\Z(G)I(H) }\bigg) \ot \bigg( \frac{\dst J}{\dst I(H)J} \bigg)
\:. $\hfill\Ph
\end{satz}\vspace{3mm}

\N{\bf Proof\,:}\quad Take $U = I (H)$, $M =I (G)$, $V=N=I(H)J$ in
Corollary \ref{pushsequ}.  Then we have isomorphisms
\vspace{2mm}

\N\makebox[14.7cm]{ \makebox[0mm]{
\begin{minipage}{17cm}\small
  \[ \frac{\dst I(H) \ot_H  J}{\dst \Imm{I(H) \ot_H I(H)J}} \hspace{2mm} \cong  \hspace{2mm}
I(H) \hspace{1mm}\ot_H \hspace{1mm}  \frac{\dst J}{\dst I(H)J}  \hspace{2mm} \cong  \hspace{2mm}
\frac{\dst I(H)}{\dst I^2(H)}  \hspace{1mm}\ot_H \hspace{1mm}  \frac{\dst J}{\dst I(H)J}  \hspace{2mm} \cong  \hspace{2mm}
\frac{\dst I(H)}{\dst I^2(H)}  \hspace{1mm}\ot \hspace{1mm}  \frac{\dst J}{\dst I(H)J} 
\]
\end{minipage}\ruled
}\rule[-11mm]{0mm}{3mm} }

\N where the first two follow from right-exactness of the tensor product and the latter is due to the fact that both $I(H)/I^2(H)$
and
$J/I(H)J$ are trivial  $H$-modules. Similarly,
\[ \frac{\dst I(G) \ot_H  J}{\dst \Imm{I(G) \ot_H I(H)J}} \hspace{2mm} \cong  \hspace{2mm}
\frac{\dst I(G)}{\dst I(G)I(H)}  \hspace{1mm}\ot \hspace{1mm}  \frac{\dst J}{\dst I(H)J \ruled} \:.
\]
Next we recall that the sequence
 \BE\label{Fox1sequ} 0 \to\: \frac{\ruleu\dst I(H)}{\dst I^2(H)} \mr{i_1} \frac{\dst I(G)}{\dst I(G) I(H)}
\mr{q_1}   \frac{\dst I(G)}{\ruled\dst \Z(G) I(H)} \:\to 0  \EE
of abelian groups is exact by the classical Fox theorem (i.e.\ Proposition \ref{Fox1} with $J=
I(H)$). Using this and the above isomorphisms diagram \REF{pushsequdia} (without the terms $K$) identifies with
the   following   diagram where in turn the bottom row is exact and induced by the top row.
  \[
\begin{matrix}
  \frac{\ruleu\dst I(H)}{\dst I^2(H)} \otz  \frac{\dst J}{\dst I(H) J}
& \mr{i_1\ot id}  & \frac{\dst I(G)}{\dst I(G) I(H)} \otz  \frac{\dst J}{\dst I(H) J}  &
\Sur{q_1 \ot id}  & \frac{\dst I(G)}{\dst \Z(G) I(H)} \otz  \frac{\dst J}{\ruled\dst I(H) J}
\cr 
  \surdown{h}  &  &  \surdown{g} &  & \| \cr
\frac{\ruleu \dst I(H)J}{\dst  I^2(H) J} & \mr{j_1}  & 
\frac{\rule{0mm}{6mm} \dst I(G)J}{\dst  I(G)I(H) J}   & \Sur{\overline{q_1\ot id}}  &  
\frac{\dst I(G)}{\dst \Z(G) I(H)} \otz  \frac{\dst J}{\ruled\dst I(H) J} 
\end{matrix}\]
  But  sequence \REF{Fox1sequ} is split for $I(G)/ \Z(G) I(H) $ is a subgroup of $\Z(G)/ \Z(G) I(H)  \linebreak
\:\cong\: \Z(G) \ot_H \Z$ which is a free abelian group since $\Z(G)$ is a free $H$-module.
Thus the top sequence in the above diagram is $R$-split short exact, whence so is the bottom 
sequence, as was to be shown.\hfill $\Box$\V

\N In the sequel we study cases of Corollary \ref{pushsequ} where $V=N=J$. Taking $U=I(H)\cap M$ it amounts to the following
description of $I(G)J/MJ$ as a group extension.\V

\begin{prop}\label{extforM} Let $M$ be a right $H$-submodule of $I(G)$. Then there is a natural commutative diagram of
$R$-modules with exact rows where
$T= {\rm Tor}_1^H(I(G)/(I(H)+M)\,,  J)$ 
and  where the top row is part of the
corresponding long exact sequence.
\vspace{2mm}

\N\makebox[14.7cm]{ \makebox[0mm]{
\begin{minipage}{17cm}\small
\BE \label{extforMdia} \begin{matrix}
T & \mr{\tau} &   \frac{\rule{0mm}{6mm}\dst I(H)   }{\dst   I(H)   \cap M  \ruled }
\hspace{1mm}  \ot_H \hspace{1mm}  J  
& \mr{i\ot id}  &  \frac{\dst I(G)}{\dst M\rule{0mm}{2.7mm}}  \hspace{1mm} \ot_H \hspace{1mm}  J   & \mr{q\ot id} &
  \frac{\dst I(G)}{\dst  I(H)+M}  
\hspace{1mm} \ot_H \hspace{1mm}  J & 
\to 0 \cr
\| & & \surdown{\mu_H} & & \surdown{\mu_G} & & \| & \cr
 T &
\mr{\mu_H  \tau} &
\frac{\rule{0mm}{6mm}\dst I(H) J  }{\dst  (I(H)   \cap M) J }  &  \mr{j} &
  \frac{\dst I_R(G) J}{\dst M_R J   \rule{0mm}{2.7mm}\ruled}&
\mr{\overline{q\ot id} }&
  \frac{\dst I(G)}{\dst  I(H)+M} 
\hspace{1mm} \ot_H \hspace{1mm}   J & \to 0 \rule[-11mm]{0mm}{3mm}
\end{matrix}\EE
\end{minipage}\ruled
}\rule[-11mm]{0mm}{3mm} }
Consequently, the bottom sequence is induced by the top one.\hfill$\Box$
\end{prop}

For $R=\Z$ and $J=I(H)$ we can improve Proposition \ref{extforM} as follows.

\begin{satz}\label{H2sequ}\quad Let $M \subset I(G)$ be a right $H$-submodule. Then there
is a natural long exact sequence
\\
\N\makebox[14.7cm]{ \makebox[0mm]{
\begin{minipage}{20cm}\small
\[ \begin{matrix}
     \frac{\ruleu\dst I^2(H)}{\dst \ruled (I(H)
\cap M)I(H)}  &  \mr{j} &  \frac{\dst  I(G)I(H)}{\dst M I(H)} &  \mr{\tilde{q}} & 
\bigg( \frac{\dst \ruleu I(G)}{\dst I(H) + M} \bigg) \ot_H \:I(H) & \to & 0 \cr
\mapup{  \mu_H\tau_1\tau_2 } & & & & & & \cr
  H_2(H,
\frac{\dst \ruleu I(G)}{\dst I(H) + M})  & \ml{\delta} & 
 H_1(H, M) & \ml{\iota_{\ast}} & H_1(H, I(H)\cap M) & \ml{} & \ldots
\end{matrix}\]
\end{minipage}\ruled
}  }
\rule{0mm}{7mm}

\N where $\T{q}(\overline{xy}) = \bar{x} \ot \bar{y}$ for $x\in I(G)$, $y\in I(H)$,  and where the dots ``$\ldots$" stand 
for extension by the long exact homology sequence induced  by the short
exact sequence $I(H) \cap M \Inj{\iota} M \auf M/I(H)\cap M$ of right $H$-modules. The other operators are defined in the
proof below.
\end{satz}

Note that the homology of $H$ here refers to $right$ coefficient modules, i.e.\ 
$H_i(H,A) = {\rm Tor}_i^H(A,\Z)$ which coincides with usual homology of the opposite module: ${\rm Tor}_i^H(A,\Z) = {\rm
Tor}_i^H(\Z,A^*)$ where $A^*$ is $A$ as an abelian group endowed with the left $H$-action $h\cdot a =ah^{-1}$, $h\in H$, $a\in A$.

We point out that the description of the {\em cokernel}\/  of $j$ given by the theorem can be used to study {\em Fox quotients}\/ 
while the description of the {\em kernel}\/ of $j$  can be used to study {\em Fox subgroups}\/, as follows.\V


%

%

%
%

\begin{kor}\label{Foxsgp}\quad Let $M \subset I(G)$ be a right $H$-submodule. Then
  \[  G \cap (1 + M I(H)) = H  \cap \Big(1 + \pi^{-1}\Imm{\mu_H\tau_1\tau_2}\Big)   \]
where $\sepi{\pi}{I^2(H)}{I^2(H)/(I(H)\cap M)I(H)}$ is the canonical projection.
\end{kor}

This means that the study of {\em generalized Fox subgroups}\/ $G \cap (1 + M I(H))$ is  
reduced to homology and induced subgroups of $H$, the latter with respect to right ideals in
$I(H)$ instead of $I(G)$.
 In particular we derive the following remarkable fact.\vspace{2mm}

\begin{kor}\label{Foxsgpfree}\quad Let $H$ be a free subgroup of $G$ and 
$M \subset I(G)$ be a right $H$-submodule. Then
 \[  G \cap (1 + M I(H)) = H  \cap (1 + (I(H) \cap M)I(H))\:.\]
\end{kor}

This is immediate from Corollary \ref{Foxsgp} noting that free groups are of cohomological dimension
1.   For specific $M$, in particular for $M = I^n(G)$, this situation will be further analyzed
in subsequent work.\V

\proofof{Theorem \ref{H2sequ}}   In Proposition \ref{extforM} take $R=\Z$ and  $J=I(H)$, and write $\tau=\tau_1$.
Now consider the following
anticommutative square consisting of connecting homomorphisms, where the vertical ones are
induced by the short exact sequence $I(H) \stackrel{\nu}{\hra} \Z(H) \Sur{\epsilon} \Z$ and the bottom one   by the sequence
$ I(H)/I(H) \cap M \hra I(G)/M \auf I(G)/I(H)+M$.
  \[\begin{matrix}
\ruled {\rm Tor}_1^H(I(G)/(I(H)+M)\,, I(H)) & \mr{\tau_1}  &  (I(H)/I(H) \cap M)\ot_H I(H)\cr
\cong\hspace{-1mm}\mapup{\tau_2}  &  &  \mapup{\tau_2\st}  \cr
{\rm Tor}_2^H(I(G)/(I(H)+M)\,, \Z) & \mr{\tau_1\st}  & 
{\rm Tor}_1^H(I(H)/I(H) \cap M,\Z) \end{matrix}\]
Here $\tau_2$ is an isomorphism since ${\rm Tor}_i^H(-,\Z(H)) = 0$ for $i\ge 1$ as $\Z(H)$ is a free $H$-module. 
Hence $\Ker{j} = \Imm{\mu_H\tau_1 \tau_2} = \Imm{\mu_H\tau_2\st \tau_1\st}$. So it
remains to exhibit $ \Ker{\mu_H\tau_2\st \tau_1\st}$. Consider the
following commutative diagram of right $H$-modules where $U=I(H)\cap M$ and $\epsilon\st(\bar{x}) =\overline{x-\epsilon(x)}$ for
$x\in
\Z(G)$ ($\epsilon\st$ is $H$-linear as for $h\in H$, $xh - \epsilon(xh) =(x-\epsilon(x))h +  \epsilon(x)(h-1)
\:\equiv\:(x-\epsilon(x))h$ mod $I(H)$).
\[\begin{matrix}
\ruled \Z(H)/U & \Inj{i'}  &  \Z(G)/M  &  \Sur{\epsilon\st}  &  I(G)/(I(H) + M) \cr
\Injup{\bar{\nu}}   &  &  \Injup{\bar{\nu}}  &  &  \| \cr
\ruleu I(H)/U & \Inj{i}  &  I(G)/M  &  \Sur{q}  &  I(G)/(I(H) + M)
\end{matrix}\]
The rows are short exact, so applying the functor ${\rm Tor}_1^H(-,\Z)$ provides the first two commutative squares of the
following diagram with exactness in the second term from the left.
\\

\N\makebox[14.7cm]{ \makebox[0mm]{
\begin{minipage}{20cm}\small
\[ \begin{matrix}
{\rm Tor}_2^H(\frac{\dst \Z(G)}{\dst M},\Z)  & \mr{\epsilon\st_{\ast}}  & 
{\rm Tor}_2^H(\frac{\dst I(G)}{\dst I(H)+M}\,,\Z) &
\mr{\tau_3}  &  {\rm Tor}_1^H(\frac{\dst \Z(H)}{\dst U},\Z)  &  \Inj{\tau_4}  & 
\bigg(\frac{\dst \Z(H)}{\dst U}\bigg)\ot_H I(H)
\cr
\mapup{\bar{\nu}_{\ast}}  &  &  \| & & \mapup{\bar{\nu}_{\ast}}& &
\mapup{\bar{\nu}\ot id}\cr 
\ruleu {\rm Tor}_2^H(\frac{\dst I(G)}{\dst M},\Z)  & \mr{q_{\ast}}  & 
{\rm Tor}_2^H(\frac{\dst I(G)}{\dst I(H)+M}\,,\Z) &
\mr{\tau_1\st}  &  {\rm Tor}_1^H(\frac{\dst I(H)}{\dst U},\Z)  &  \mr{\tau_2\st}  & 
\bigg(\frac{\dst I(H)}{\dst U}\bigg)\ot_H I(H)
\end{matrix}\]
\end{minipage}\ruled
}\rule[-5mm]{0mm}{3mm} }
\rule{0mm}{7mm} \vspace{-1mm}

\N The operators $\tau_2\st$ and $\tau_4$ are injective as they are part of the corresponding long exact sequences induced by the
short exact sequence
$I(H) \Inj{\nu} \Z(H) \Sur{\epsilon} \Z$, and as ${\rm Tor}_1^H(-,\Z(H))=0$. 
Now under the isomorphisms $(\Z(H)/U) \ot_H I(H) \:\cong\:  \Z(H) \ot_H I(H)/\Imm{U \ot_H
I(H)}  \:\cong\: I(H)/UI(H)$ the map $(\bar{\nu} \ot id)$ corresponds to $\mu_H$, so that 
$\Ker{\mu_H\tau_2\st \tau_1\st} = $ Ker$((\bar{\nu}\ot id)\tau_2\st \tau_1\st) =  \Ker{\tau_4\tau_3}
=  \Ker{\tau_3} = \Imm{\epsilon\st_{\ast}}$.
To compute \Ker{\epsilon\st} we consider the canonical extension to the left of the top row of the foreging diagram, and we use the
connecting homomorphisms induced by the short exact sequences $U\hra \Z(H) \auf \Z(H)/U$ and $M\hra \Z(G) \auf \Z(G)/M$ which are
isomorphisms as $\Z(H)$ and $\Z(G)$ are free $H$-modules.
  \[ \begin{matrix}\ruleu \ruled
{\rm Tor}_1^H(U,\Z)  &  \mr{\iota_{\ast}} &  {\rm Tor}_1^H(M,\Z) &  & \cr
\isoup{\tau_5}  &  & \isoup{\tau_6} & &  \cr
\ruleu {\rm Tor}_2^H(\frac{\ruleu\dst \Z(H)}{\dst U},\Z)  & \mr{i\st_{\ast}} & 
{\rm Tor}_2^H(\frac{\dst \Z(G)}{\dst M},\Z)  & \mr{\epsilon\st_{\ast}}  & 
{\rm Tor}_2^H(\frac{\dst I(G)}{\dst I(H)+M}\,,\Z) \:.
\end{matrix}\]
So putting $\delta =  \epsilon\st_{\ast} \tau_6^{-1}$ the theorem is proved. \hfbox\V

We still give another application of Proposition \ref{extforM} which is motivated by the computation of the induced subgroup $G\cap
\Big(1 + I(H)I(G)I(H) + I([H,G])I(H)
\Big)$   by Tahara, Vermani and Razdan in \cite{Razdan}; we show that their result can be recovered and generalized to
arbitrary $R$ by the computation of a related quotient of $I_R(G)$, as follows.\V

\begin{satz}\label{KVgen} Suppose that $H$ is normal in $G$. Then there is a non-natural isomorphism of $R$-modules 
  \[  \frac{\dst I_R(G) J}{\dst I_R(H)I_R(G) J +  I_R([H,G]) R(G)J}   \hspace{2mm}\cong\hspace{2mm}  
\Big(\frac{\dst H}{\dst [H,G]}  \hspace{1mm}\ot\hspace{1mm} \frac{\dst
J}{\dst I(H)J}\Big)\hspace{2mm}\oplus\hspace{2mm}  \Big(I(G/H) \hspace{1mm}\ot\hspace{1mm} \frac{\dst
J}{\dst I(H)J}\Big)
\] 
as well as identities  
  \begin{eqnarray*}
\lefteqn{  G\cap \Big(1 + I_R(H)I_R(G)I_R(H) + I_R([H,G])I_R(H) \Big) \hspace{50mm} \phantom{i}}\\ 
  \hspace{20mm} &=& G\cap \Big(1 +  I_R(H)I_R(G) I_R(H) +  I_R([H,G])R(G)I_R(H) \Big) \\
  &=&  H \cap \Big(1 +  I_R([H,G])I_R(H) + I_R^3(H)\Big) \:. \\
\end{eqnarray*}
\end{satz}\vspace{-3mm}

  We point out that the latter term can be easily made explicit by using our formula in \cite{D3F2} for the second relative
dimension subgroup
$G\cap (1+I_R(K)I_R(G)+I_R^3(G))$ for all $K\le G$ and $R$ (reproved differently for $R=\Z$ in \cite{Razdan}); as the result is
rather involved for arbitrary $R$ we refrain from stating it here.\V

\proof In Proposition \ref{extforM}, take 
$M= I(H)I(G) + I([H,G])\Z(G)$ and 
$V=N=J$. Let $N$ a normal subgroup of $G$. Then there is a canonical isomorphism of right $G$-modules 
  \begin{equation}\label{psiNdef} 
  \Psi_N\:\colon\: I(G)/I(N) \Z(G)  \hspace{2mm}\to\hspace{2mm} I(G/N)\:,\quad \overline{a-1} \mapsto aN-1\:.
   \end{equation}
In particular, the map $\Psi_{[H,G]}$ induces another  isomorphism of right $G$-modules 
$\gamma\,\colon\: I(G)/M    \:\stackrel{\cong}{\lra}\: I(G/[H,G])/ I(H/[H,G]) I(G/[H,G])$. Furthermore, $N$ gives rise to the
following  exact sequence of left $G/N$-modules where $D(\bar{a}) = \overline{a-1}$ for $a\in H$, see Theorem VI.6.3 in
\cite{Hi-St}.
  \BE\label{Foxcalculsequ}
0 \hspace{1mm}\to \hspace{2mm}N\ab \hspace{2mm}\mr{D}\hspace{2mm} \frac{\dst I(G)}{\dst I(N)I(G)}
\hspace{2mm} \mr{\overline{I(q\st)}} \hspace{2mm} I(G/N)\hspace{2mm}
\to \hspace{1mm} 0 \:.\EE 
In particular, taking $N=G$ one gets the wellknown isomorphism of groups
\BE\label{IG/I2G} \phi_G\,\colon\, G/G_2 \cong I(G)/I^2(G)\:, \quad \phi_G(aG_2) = a-1+I^2(G) \quad \mbox{for $a\in G$.}\EE
Taking $N=[H,G]$ we obtain the following commutative diagram with exact  rows where $U=I(H) \cap M$.
  \[\begin{matrix}
\ruleu\ruled 0  & \to & \frac{\dst I(H)}{\dst \ruled\rule{0mm}{3.5mm} U}  & \mr{i}  & \frac{\dst I(G)}{\dst \rule{0mm}{3.5mm}M}  & 
\mr{q}  &  \frac{\dst I(G)}{\dst
\rule{0mm}{3.5mm} I(H)+M }  & \to & 0 
\cr
  &  & \mapdown{\tilde{\gamma}}  &  &  \isodown{\gamma}  &  &  \isodown{\Psi_H} \cr
\ruleu\ruled 1  & \to & H/[H,G]  &  \mr{D}  &  \frac{\dst \ruleu\rule{0mm}{-4mm}I(G/[H,G])}{\dst \rule{0mm}{3.5mm} I(H/[H,G])
I(G/[H,G])}  & 
\mr{\overline{I(q\st)}}  &  I(G/H)  & \to & 0 \end{matrix}\]
Here $\tilde{\gamma}$ is induced by $\gamma$ and hence is an isomorphism, too. The bottom sequence is \Z-split as $I(G/H)$ is a free
\Z-module; hence so is the top sequence.  Moreover,
the right $H$-action on each of its terms is trivial; for $I(H)/U$ and $I(G)/M+I(H)$ this follows from the inclusions $I^2(H)
\subset U$ and $\Z(G)I(H) = I(H)\Z(G) \subset M+I(H)$, resp., and for $I(G)/M$ by right $G$-linearity of $\gamma$ and the fact that
$H/[H,G]$ is central in $G/[H,G]$, whence $I(H/[H,G]) I(G/[H,G]) = I(G/[H,G]) I(H/[H,G]) $. Using these facts diagram
\REF{extforMdia} (without the terms $T$) identifies with the following diagram.
\\
\N\makebox[14.7cm]{ \makebox[0mm]{
\begin{minipage}{20cm}\small
  \[\begin{matrix}
\ruleu\ruled 0  & \to & \frac{\dst I(H)}{\dst \ruled\rule{0mm}{3.5mm}U} \ot \frac{\dst J}{\dst \rule{0mm}{3.5mm}I(H)J}  & \mr{i\ot id} 
&  \frac{\dst I(G)}{\dst \rule{0mm}{3.5mm}M}  \ot
\frac{\dst J}{\dst \rule{0mm}{3.5mm}I(H)J}  & 
\mr{q\ot id}  &  I(G/H)  \ot \frac{\dst J}{\dst \rule{0mm}{3.5mm}I(H)J}  &
\to & 0 
\cr
  &  & \surdown{\bar{\mu}_H}  &  &  \surdown{\bar{\mu}_G}  &  &  \| \cr
\ruleu\ruled  0 & \to & \frac{\dst \ruleu I(H)J}{\dst I([H,G])J + I^2(H)J}  &  \mr{j}  &  \frac{\dst I(G) J}{\dst
\rule{0mm}{3.5mm} I(H)I(G) J + I([H,G])
\Z(G)J}  &  \mr{\overline{q\ot id}}  &  I(G/H)  \ot \frac{\dst J}{\dst \rule{0mm}{3.5mm}I(H)J}  & \to & 0 
 \end{matrix}\]
\end{minipage}\ruled
}\rule[-24mm]{0mm}{3mm} }

\N By the foregoing the top row is $R$-split exact, whence so is the bottom row by Proposition \ref{extforM} which implies the
desired isomorphism. Now take $J=I_R(H)$. We claim that the following relations hold.
  \begin{eqnarray}
\lefteqn{G\cap \Big(1 + I_R(H)I_R(G)I_R(H) + I_R([H,G])I_R(H) \Big)} \nonumber \\  &\subset&     G\cap \Big(1 +  I_R(H)I_R(G)
I_R(H) + I_R([H,G])
R(G)I_R(H) \Big)
\label{Razdanequ1}\\
  &=&  H\cap \Big(1 +  I_R(H)I_R(G) I_R(H) + I_R([H,G]) R(G)I_R(H) \Big)  \label{Razdanequ2}\\
  &=&  H \cap \Big(1 +  I_R([H,G])I_R(H) + I_R^3(H)\Big)  \label{Razdanequ3}
\end{eqnarray}
Indeed, \REF{Razdanequ1} being plain, \REF{Razdanequ2} follows from the  relation
$G\cap (1+I_R(G)I_R(H))$ $ = H\cap (1+I^2_R(H))$ (cf.\ Corollary \ref{GcapIJprop} below), which also implies \REF{Razdanequ3}
by injectivity of the map 
$j$ in the above diagram.  As the last term is contained in the first the  identities of induced subgroups given in the assertion
are proved.\hfbox\V


\section{Fox polynomial groups}

Polynomial groups $P_{n,R}(G) = I_R(G)/I_R^{n+1}(G)$ were introduced by Passi in \cite{Pa68} (see also \cite{Pa}), along with a
notion of  polynomial maps  from groups to $R$-modules such that the map $p_{n,R}\,\colon\,G \to P_{n,R}(G)$, $p_{n,R}(a) =
a-1+I_R^{n+1}(G)$ is universal polynomial of degree $\le n$. {\em Relative}\/ polynomial groups
$P_{n,R}(G,N) = I_R(G)/(I_R(N)I_R(G) + I_R^{n+1}(G))$ for normal subgroups
$N$ were formally introduced in \cite{PolProp} but had been implicitely studied in the literature before. Indeed, these
constructions proved to be very useful in the study of dimension subgroups $D_{n,R}(G) :\,= G\cap (1+I_R^n(G)) =
\Ker{p_{n-1,R}}$ as well as of augmentation quotients
$Q_{n,R}(G)
:\,= I_R^{n}(G)/I_R^{n+1}(G) = \Imm{(G\ab)^{\ot n}\ot R \to P_{n,R}(G)}$, see \cite{Pa} and \cite{Q3}; moreover, they are
used in
\cite{PolProp} to study extensions of torsionfree nilpotent groups and in \cite{GoG} to determine the Schur multiplier of $2$-step
nilpotent groups.  We here extend the  approach of augmentation quotients via polynomial groups to Fox quotients $Q_{n,R}(G,H) =
I_R^{n-1}(G)I_R(H)/I_R^{n}(G)I_R(H)$: first consider {\em relative Fox polynomial groups}\/ 
   \[ P_{n,R}(G,N;H) = R(G)I_R(H) \Big/\Big( R(G)I_R(N)I_R(H) + I_R^{n}(G)I_R(H) \Big) \hspace{1mm}, \]
 \[ \bar{P}_{n,R}(G,N;H) = I_R(G)I_R(H) \Big/\Big( R(G)I_R(N)I_R(H) + I_R^{n}(G)I_R(H) \Big) \hspace{1mm}; \]
note that $P_{n,R}(G,N) =  P_{n,R}(G,N;G)$ and  $P_{n,R}(G,N)^2 =  \bar{P}_{n,R}(G,N;G)$. In a second step
 determine $Q_{n,R}(G,H)$ as a subgroup of $\bar{P}_{n,R}(G,1\!\!1;H)$. We here carry out this program for $n=2$ and all $R$
and for $n=3$ and $R=\Z$.

Actually we consider a more general version of $\bar{P}_{n,R}(G,K;H)$, replacing $I_R(H)$ by $J$ and the augmantation powers
$I_R^{n}(G)$ by the filtration terms $\IRN{n}$ induced by an N-series $\cal G$ of $G$. This degree of generality is necessary in
many contexts, such as polynomial cohomology \cite{PolProp} or (ordinary) augmentation and Fox quotients and subgroups of
semidirect products \cite{Ta}, \cite{Kh12}, and even of arbitrary groups:   we show in subsequent work that under
suitable conditions  the classical Fox subgroup $G\cap \big(1+ I^n(G)I(H)\big)$ equals $H\cap \big(1+ I^n_{Z\!\!\!Z,\cal H}(H)I(H)\big)$ for
an appropriate N-series ${\cal H}$ of $H$. 

Let us recall the necessary definitions. An {\em N-series}\/ ${\cal G}$  of $G$
is a descending chain of subgroups
  \[   G= G_{(1)} \supset G_{(2)}\supset G_{(3)}\supset \ldots  \]
such that $[ G_{(i)}, G_{(j)}] \subset G_{(i+j)}$ for $i,j\ge 1$.
A given  N-series ${\cal G}$ induces a descending chain of two-sided ideals of $R(G)$
  \[  R(G) = I^0_{R,\cal G}(G) \supset I^1_{R,\cal G}(G) \supset I^2_{R,\cal G}(G) \supset \ldots  \]
by defining $I^n_{R,\cal G}(G)$  (for $n\ge 1$) to be the $R$-submodule of $R(G)$ generated by the elements
  \[  (a_1 -1) \cdots (a_r -1)\,, \quad \mbox{$r\ge 1$, $a_i \in N_{k_i}$, such that $k_1+ \ldots +
k_r  \ge n$.}  \] 
For example, the following N-series frequently appear in the literature: 
\begin{itemize}

\item the lower central series, denoted by $\gamma=\gamma_G$, where $I^n_{R,\gamma}(G) = I^n_R(G)$;

\item the series $\sqrt{\gamma}$ defined by $G_{(i)} = \sqrt{G_i}$, the isolator of $G_i$; here $I^n_{Z\!\!\!Z,\sqrt{\gamma}}(G)$
is the isolator of $I^n(G)$, see \cite{PolProp};

\item if $G$ itself is a subgroup of a group $\Gamma$, the N-series defined by  $G_{(i)} = G\cap \Gamma_i$ or $G_{(i)} =
[G_{(i-1)},\Gamma]$ if $G$ is normal in $\Gamma$, see \cite{Gu}, \cite{Ta}.  

\end{itemize}

In the sequel we add a supercript $\cal G$ to the terms defined above (writing $P_{n,R}^{\cal G}(G)$, $P_{n,R}^{\cal G}(G,N;H)$
etc.) when we replace $I_R^k(G)$, $k=n-1,n,n+1$, by $I_{R,\cal G}^k(G)$ in their definition; and the super- or
subscripts $R,\cal G$ are frequently suppressed from the notation when $R=\Z$ or  ${\cal G} =\gamma$, resp.

Actually the pushout lemma and its corollary in section 1  provide the key tools for studying these notions:
it turns out that the isomorphism
$Q_2^{\cal G}(G) = P_2^{\cal G}(G)^2 +p_2(G_{(2)})\:\cong\: {\rm U}_2{\rm L}^{\cal G}(G)$ and its analogon for $P_3^{\cal
G}(G)^2+p_2(G_{(3)})$ in
\cite{Q3} formally generalize to the Fox case once written in form of a pushout, see Theorems \ref{Fox2IH}
and \ref{Fox3push} below. At the same time the method provides a strictly dual
approach to Fox subgroups which is pursued in
\cite{D3F2} to determine the second  relative Fox  subgroup over any ring of coefficients.\V

We start by an easy case, noting that there is a canonical approximation of the Fox quotient
$Q_n^{\calG}(G,H) =I^{n-1}_{\calG}(G)I(H)/I^{n}_{\calG}(G)I(H)$, namely 
by means of the natural surjective homomorphism
  \BE\label{zetadef} \sepi{\zeta_n^{\calG}}{Q_{n-1}^{\calG}(G)\ot H\ab \hspace{2mm}}{\hspace{2mm} Q_n^{\calG}(G,H)} \EE
where
$\zeta_n( (x+I^n_{\calG}(G)) \ot hH_2 ) = x(h-1) + I^n_{\calG}(G)I(H)$ for $x \in I^{n-1}_{\calG}(G)$, 
$h\in H$. This map rarely is an isomorphism, but it is at least in the following case. 

Let $\bar{\cal G}$ be the N-series of $G/N$ given by the image of $\cal G$ under the canonical projection.

\begin{prop}\label{Hfreequot} If $H$ is a free group then there are isomorphisms
  \[ \bar{P}_{n,R}^{\calG}(G,N;H) \hspace{2mm}\cong\hspace{2mm} P_{n-1,R}^{\bar{\cal G}}(G/N) \ot H\ab   \]
  \[  Q_{n,R}^{\calG}(G,H) \hspace{2mm}\cong\hspace{2mm} Q_{n-1,R}^{\calG}(G) \ot H\ab  \]
 the first of which is non natural while the second one is given by $\zeta_n^{\calG}$.
\end{prop}

\proof  The first isomorphism is immediate from Corollary \ref{Hfree} by taking $M=\Z(G)I(N) + I^n(G)$ and using the isomorphism $\Psi_N$ in \REF{psiNdef}. Then the composite map
$Q_{n-1,R}^{\calG}(G) \ot H\ab \Sur{\zeta_n^{\cal G}} Q_n^{\calG}(G,H) \hra
\bar{P}_{n,R}^{\calG}(G,\{1\};H) \hspace{2mm}\cong\hspace{2mm} 
P_{n-1,R}^{ \cal G }(G) \ot H\ab$ equals $i\ot id$ with $i\,\colon\,Q_{n-1,R}^{\calG}(G) \hra P_{n-1,R}^{ \cal G }(G)$. But 
$i\ot id$ is injective as $H\ab$ is a free \Z-module, hence so is $\zeta_n^{\cal G}$.\hfbox\V

For a refinement of this fact for free {\em nilpotent}\/ groups $H$ see  section 4 below.\V

To study the structure of $I_R(G) J/(R(G)I_R(N)J + \IRN{n} J)$ for arbitrary $H$ we use Lemma \ref{pushlem} for   $M=I(N)\Z(G) +
I^n_{\calG}(G)$ and
$U$ being some right
$H$-submodule of $I(H)$ and $M$ containing
$I^n(H)$; one may take $U$ equal to $I^n(H)$, cf.\ the proofs of corollaries \ref{Hfreeabcyc}, \ref{Fox2IH} and Theorem
\ref{Fox3push}, or equal to $I(H\cap N)\Z(H) + I^n(H)$ where $I(H)/U \:\cong\:P_{n-1}(HN/N)$, or equal to 
$I(H)
\cap M$, cf.\ Proposition \ref{allgsequ}. Now consider the following commutative square.
  \BE\label{diasatz}
 \begin{matrix}
\bigg(  \frac{\ruleu\dst I(H)}{\dst U} \bigg) \ot_H \bigg( \frac{\dst J}{\dst
I^{n-1}(H) J}  \bigg)  &   \mr{i\ot id}  & P_{n-1}^{\bar{\cal G}}(G/N)
\ot_H \bigg( \frac{\dst J}{\ruled\dst I^{n-1}(H) J} \bigg)   \cr 
\surdown{\mu_H}  & & \surdown{\mu_G} \cr 
 \frac{\ruleu\dst I_R(H) J}{\dst U_R J}  & \mr{j} 
& \frac{\dst I_R(G) J}{\dst R(G)I_R(N)J + \IRN{n} J \ruled} \end{matrix}\EE\vspace{1.5mm}

\begin{prop}\label{satz}\quad Diagram \REF{diasatz} is a pushout square of abelian
groups; in other words, the map $\mu_G$ induces a natural isomorphism
\[  \frac{\dst I_R(G) J}{\dst R(G)I_R(N)J + \IRN{n} J } \:\cong\:  P_{n-1}^{\bar{\cal G}}(G/N) \ot_H \bigg( \frac{\dst J}{\dst I^{n-1}(H)
J}\bigg) 
\bigg/(i\ot id){\rm Ker}({\mu_H}) \:.\]
This also implies the relation\V

  \N\Ph\hfill$  \Ker{j} = \mu_H\, \Ker{i \ot id}\:.$
\hfill\Ph
\end{prop}\vspace{3mm}

\N{\bf Proof\,:}\quad This is  a special case of Lemma \ref{pushlem}: Take  $V=N=J$ in \REF{dia1}. Then right-exactness of
the tensor product provides isomorphisms
  \BE\label{tensquot}   \frac{\ruleu\dst I(G)\ot_H J}{\dst \mbox{Im}(M \ot_H J)} \:\cong\: 
 \bigg( \frac{\dst I(G) }{\dst M}  \bigg) \ot_H J  \:\stackrel{}{\cong}\: 
\bigg(  \frac{\dst I(G) }{\dst M}  \bigg) \ot_H   \bigg( \frac{\dst J}{\ruled\dst
I^{n-1}(H)J} \bigg) \:, \EE similarly for $G$ being replaced by $H$.
Using this it is easy to identify diagram \REF{diasatz} with diagram \REF{dia1}, so Lemma
\ref{pushlem} gives the result. \hfill $\Box$\V


\N Proposition \ref{extforM} here takes the following form.\V

\begin{kor}\label{allgsequ}\quad  For $M=  I(N)\Z(G) + I^n_{\calG}(G)$ there is a natural commutative diagram with exact rows 
where we abbreviate $U=I(H)\cap M$, $J_{n-1} = J/I^{n-1}(H)J$, $T= {\rm Tor}_1^H(I(G)/(I(H)+M)\,, J_{n-1})$,  $\bar{H}
=HN/N$ and  where the top row is part of the
corresponding long exact sequence.
\\
\N\makebox[14.7cm]{ \makebox[0mm]{
\begin{minipage}{20cm}\small
\[ \begin{matrix}
T & \mr{\tau} &   \frac{\rule{0mm}{6mm}\dst I(H)   }{\dst   U  }
\hspace{1mm}  \ot_H \hspace{1mm}  J_{n-1}  
& \mr{i\ot id}  &  P_{n-1}^{\bar{\cal G}}(G/N)  \hspace{1mm} \ot_H \hspace{1mm}
  J_{n-1}   & \mr{q\ot id} &
  \frac{\dst P_{n-1}^{\bar{\cal G}}(G/N)}{\dst \langle p_{n-1}(\bar{H}) \rangle \rule{0mm}{3.5mm}}  
\hspace{1mm} \ot_H \hspace{1mm}   J_{n-1} & 
\to 0 \cr
\| & &  \surdown{\mu_H} & &  \surdown{\mu_G} & & \| & \cr
 T &
\mr{\mu_H  \tau} &
\frac{\rule{0mm}{6mm}\dst I(H) J  }{\dst  (I(H)   \cap M) J }  &  \mr{j} &
  \frac{\dst I_R(G) J}{\dst R(G)I_R(N) J + \IRN{n}J  \ruled}&
\mr{\overline{q\ot id} }&
  \frac{\dst P_{n-1}^{\bar{\cal G}}(G/N)}{\dst \langle p_{n-1}(\bar{H}) \rangle \rule{0mm}{3.5mm}} 
\hspace{1mm} \ot_H \hspace{1mm}   J_{n-1} & \to 0 \rule[-11mm]{0mm}{3mm}
\end{matrix}\]
\end{minipage}\ruled
}\rule[-11mm]{0mm}{3mm} }
Consequently, the bottom sequence is induced by the top one.\hfill$\Box$
\end{kor}

\begin{kor} \label{Hfreeabcyc} (i) There is a natural exact sequence
\BE\label{H2Hsequ} H_2(H,R) \mr{\tilde{\lambda}}  P_{n-1}^{\bar{\cal G}}(G/N) \ot_H P_{n-1,R}(H)   \mr{\mu_G} 
\bar{P}_{n,R}^{\calG}(G,N;H) \to 0 \:.\EE 

\N(ii) If $H$ is abelian and $R=\Z$,   sequence \REF{H2Hsequ} becomes
\[  H\sm H  \mr{\lambda\st} P_{n-1}^{\bar{\cal G}}(G/N) \ot_H P_{n-1}(H)  \mr{\mu_G}  \bar{P}_n(G,N;H) \to 0\]
where $H\sm H = H\ot H/\langle \{h\ot h\,|\,h\in H \}\rangle$ and $\lambda\st(h\sm h\st) = p_{n-1}(h) \ot p_{n-1}(h\st) -
p_{n-1}(h\st) \ot p_{n-1}(h)$.\V

\N(iii) Let $H=\langle t \rangle $ be cyclic of order $m$. Recall the norm element $N(t) = 1 + t + \ldots +t^{m-1}$ and let
$\bar{N}(t)$ be its image in $\Z(G/N)$. Then there is an isomorphism
  \[ \gamma\:\colon\hspace{2mm} \bar{P}_n(G,N;H) \hspace{2mm}\cong\hspace{2mm} P_{n-1}^{\bar{\cal G}}(G/N)  \Big/
P_{n-1}^{\bar{\cal G}}(G/N) \bar{N}(t)  \] which sends the coset of $xy$, $x\in I(G)$ and $y\in I(H)$, to the coset of $xy\st$
where $y\st \in \Z(H)$ such that $y=y\st(t-1)$. 
\end{kor}

\proof  Defining
$\T{\lambda} = (i\ot id)\lambda
$ assertion (i)  follows from Lemma \ref{pushlem}. Then part (ii) is an immediate consequence of the wellknown
computation 
$H_2(H)
\:\cong\:H
\sm H$ for abelian groups
$H$. Now let $H=\langle t \rangle $ be  of order $m$.  Then there is an $H$-linear isomorphism $\Z(H)/ \Z(H) N(t) \:\mr{\cong}\:
I(H)$ sending $\bar{x}$ to $x(t-1)$ for $x\in \Z(H)$. As here $H\sm H =0$ one gets isomorphisms
  \begin{eqnarray*} \bar{P}_n(G,N;H) &\cong& P_{n-1}^{\bar{\cal G}}(G/N) \ot_H P_{n-1}(H)   \\
&\cong & 
P_{n-1}^{\bar{\cal G}}(G/N) \ot_H I(H)   \\
&\cong& P_{n-1}^{\bar{\cal G}}(G/N) \ot_H  \Big( \Z(H)/ \Z(H)N(t)   \Big)  \\
&\cong &
 \Big( P_{n-1}^{\bar{\cal G}}(G/N)/P_{n-1}^{\bar{\cal G}}(G/N)\cdot\Z(H) N(t) \Big)  \ot_H  \Z(H) \\
&\cong & P_{n-1}^{\bar{\cal G}}(G/N) \Big/ P_{n-1}^{\bar{\cal G}}(G/N) \bar{N}(t)  
\end{eqnarray*}
where the dot $\cdot$ denotes the canonical right $H$-module structure of $P_{n-1}^{\bar{\cal G}}(G/N)$.\hfbox\vspace{1mm}

  We remark that parts (ii) and (iii) of the preceding corollary can easily be generalized to arbitrary coefficient rings $R$ by
using part (i) and the universal coefficient sequence, along with the explicit generators of \,Tor$_1^{Z \!\!\! Z}(H\ab,R)$
provided in
\cite{ML} V.6. We leave it to the interested reader to write out the details.\vspace{3mm}


\section{The first two generalized Fox quotients}

We first quote the following elementary fact which was first proved by Whitcomb for $J=I(H)$ in \cite{Whitcomb}.\V

\begin{prop}\label{Fox1}\quad The first generalized Fox quotient is
reduced to the coinvariants of $J$ by the natural isomorphism\V
  
  \Ph\hfill$  \frac{\dst R(G)J}{\dst I_R(G)J} \: \cong\: \frac{\dst J}{\dst I(H)J}\:. $\hfill
$\Box$ 
\end{prop}\vspace{1.5mm}

  This may be well-known but was  reproved anyway together with Lemma \ref{pushlem}.

\begin{kor}\label{GcapIJprop} If $J\subset I_R(H)$ one has the relation
\[G \cap (1 + I_R(G)J) = H \cap (1 + I_R(H)J)\,.\]
\end{kor}

This was proved for $R=\Z$ in \cite{Ka-Ve}.\V

\proof We know from \cite{Unico} that
\begin{equation}\label{GcapIJ}
G \cap (1 + I_R(G)J) = H \cap (1 + I_R(H)J)\,.
\end{equation}
Now let $T$ be a left transversal for $H$ in $G$ containing $1$. Then $R(G) = \bigoplus_{t\in T}t.R(H)$ and $R(G)J = \bigoplus_{t\in T}t.J \subset 
\bigoplus_{t\in T}t.R(H)$. As $I_R(H)$ lies in the summand $1.R(H)$ we see that $I_R(H)\cap R(G)J =J$, whence
\begin{equation}\label{HcapR(G)J}
H \cap (1 + R(G)J) = H \cap (1 + J)\,.
\end{equation}
The assertion now follows combining \REF{GcapIJ}, \REF{HcapR(G)J} and Proposition \ref{Fox1}.\hfill$\Box$\V


  Now we turn to the second generalized Fox quotient. Let $K$ be some subgroup of $G$. Take $n=2$ and $N=KG_{(2)}$   in Corollary
\ref{allgsequ}, noting that then
$M= \Z(G)I(KG_{(2)})  + I^2_{\calG}(G) = I(K) + I^2_{\calG}(G) = I(KG_{(2)}) + I^2(G) $.  All tensor products over $H$
reduce to tensor products over the integers as all terms involved are trivial $H$-modules. Another reduction comes from the  
natural isomorphism $\phi_G $ (cf.\ \REF{IG/I2G}) which induces a canonical isomorphism $P_1^{\bar{\cal G}}(G/N) =
I(G)/(I(KG_{(2)}) + I^2(G))
\:\cong\:G/KG_{(2)}$; it  implies that $I(H)\cap M = I(H\cap KG_{(2)}) + I^2(H)$ since 
\begin{eqnarray*} \frac{I(H)\cap (I(KG_{(2)}) + I^2(G))}{I^2(H)}  &=&  \mbox{Ker} \Big( I(H)/I^2(H) \mr{\bar{i}} I(G)/I(KG_{(2)}) +
I^2(G)
\Big)  
\\
  & = & \phi_H  \Ker{H/H_2 \to G/KG_{(2)}} \:.
\end{eqnarray*}
Hence $I(H)/I(H)\cap M\:\cong\: HKG_{(2)}/KG_{(2)}$. Under these identifications the diagram in Corollary \ref{allgsequ}
(except from the Tor-terms) looks as follows where $J_H = J/I(H)J$, $T\st=
{\rm Tor}_1^{\Z}(G/HKG_{(2)}\,,J_H) $ and where the top row is part of the corresponding
6-term exact sequence  with $\iota$ being
induced by the inclusion $H \hra G$ and $\pi$ being the canonical quotient map. Moreover, $\mu_H\st(\bar{h} \ot\bar{x}) =
\overline{(h-1)x}$ and $\mu_G\st(\bar{g} \ot\bar{x}) =
\overline{(g-1)x}$ for $h\in H$, $g\in G$ and $x\in J$.
\\
\N\makebox[14.7cm]{ \makebox[0mm]{
\begin{minipage}{16cm} \small
\BE\label{Fox2dia}\begin{matrix}
\hspace*{2mm} T\st & \mr{\tau}  & 
\Big( \frac{\ruleu\dst HKG_{(2)}}{\dst \ruled KG_{(2)}} \Big) \otz J_H
& \mr{\iota\ot id}  & \Big( \frac{\dst G}{\dst KG_{(2)}} \Big)\otz  J_H
& \Sur{\pi \ot id}  & \Big( \frac{\dst G}{\dst HKG_{(2)}}\Big) \otz  J_H \cr 
\|  &  & \surdown{\mu_H\st}  &  &  \surdown{\mu_G\st} &  & \| \cr
T\st & \mr{\mu_H\st \tau}  & 
\frac{\ruleu \dst I_R(H)J}{\dst I_R(H \cap KG_{(2)}) J+ I^2_{R}(H) J} & \mr{j }  & 
\frac{\ruleu \dst I_R(G)J}{\dst I_R(K) J + I_{R,{\cal G}}^2(G) J}
  & \Sur{\overline{\pi\ot id}}  &  \Big( \frac{\dst G}{\dst HKG_{(2)}}\Big)
\otz J_H
 \end{matrix}\EE
\end{minipage}$\ruled$
}}


\N\rule{0mm}{11mm} In view of the above remarks Corollary \ref{allgsequ} amounts to the following functorial description of the
$R$-module
  $I_R(G)J/(I_R(K)J + I_{R,{\cal G}}^2(G)J)$.\V

\begin{satz}\label{Fox2}\quad   The middle square of diagram \REF{Fox2dia} is a pushout of $R$-modules; in particular, 
$\Ker{\mu_G\st} = (\iota \ot id)\Ker{\mu_H\st}$. 
  Moreover, the rows of the  diagram above are exact; in particular,  $\Ker{j }$ is generated by the cosets $\overline{(h -1) x}$
where
$h
\in H$, 
$x \in J$  such that there is some $k \in
\Z$ for which $h \in KG_{(2)}G^k$ and  $kx \in I(H)J$. Furthermore,
omitting the 
terms $T\st$, the bottom sequence is induced by the top
sequence via the map $\mu_H\st$.\hfill $\Box$
\end{satz}\vspace{3mm}

The  description of $\Ker{j }$ given in the theorem is obtained by using the canonical generators of the
torsion product, see \cite{ML} V.6.\vspace{3mm}

\begin{kor}\label{Fox2split}\quad If the sequence 
  \[   1 \to  HKG_{(2)}/KG_{(2)} \to G/KG_{(2)} \to G/HKG_{(2)} \to 1 \]
of abelian groups splits then there is an isomorphism of $R$-modules\\
  \Ph\hfill$ \frac{\ruleu \dst I_R(G) J}{\ruled \dst I_R(K)J + I_{R,{\cal G}}^2(G) J}  \:\cong\:  
\frac{\dst I_R(H)J }{\dst      I_R(H
\cap KG_{(2)}) J + I_R^2(H)J}   \:\oplus\: \Big( \frac{\dst G}{\dst HKG_{(2)}}\Big)
\otz \Big(  \frac{\dst J}{\ruled\dst I(H) J}\Big)  \:.$ 
\end{kor}\vspace{1.5mm}

In fact, here the top sequence in diagram \REF{Fox2dia} is split short exact whence so is the bottom sequence (omitting the terms
$T\st$).\hfbox\V

\N This last result in particular applies  when $K$ is normal, $\calG=\gamma$  and 
$HK/K$ is a semidirect factor of $G/K$  (normal or not). 
Thus we recover the corresponding results of Khambadkone and Karan-Vermani for $J=I(H)$  
in \cite{Kh13}, \cite{KV30} and for $J=I^m(H)$ in
\cite{Kh14}, \cite{KV29}.\V

\N Another application of Theorem \ref{Fox2} is the following intersection theorem.\V

\begin{kor}\label{intersec} For any subgroup $K$ of $G$ one has the identity
  \[ I_R(H)J \hspace{2mm} \cap  \hspace{2mm} \Big(I_R(K)J + I_{R,{\cal G}}^2(G) J \Big)  =  I_R(H \cap KG_{(2)}) J  \hspace{2mm} +
\hspace{2mm}  I_R^2(H)J  
\hspace{2mm} + \hspace{2mm}  U \] where the subgroup $U$ of $I_R(H)J $ is generated by the elements $(h-1)x$ where $h\in H$, $x\in
J$ for which there exists $k\in
\Z$ such that $h\in KG_{(2)}G^k$ and $kx\in I(H)J$.\linebreak
\hfbox\end{kor}

  In particular, the groups $I^n(H) \cap \Big(I(K) I^{n-1}(H) + I^2(G)I^{n-1}(H) \Big)$ can be considered as being known for $n\le
4$, using the computation of $I^2(H)/I^3(H)$ in \cite{Ba-GrII} or \cite{Q3} and of $I^3(H)/I^4(H)$ in \cite{Q3}.

In the case $R=\Z$ and $J = I(H)$ we can improve Theorem \ref{Fox2}, as follows. For a group $K$ with N-series ${\cal K}$ there is
a canonical  homomorphism
  \[ c_2^{K,\cal K}\:\colon\hspace{2mm} (K/K_{(2)}) \sm (K/K_{(2)}) \to  K_{(2)}/K_{(3)}\:,\quad  
 c_2^{K,\cal K}((aK_{(2)}) \sm (bK_{(2)})) = [a,b]K_{(3)} \]
for $a,b\in K$. We note $c_2^{G,\cal G} =c_2^{\cal G}$ and $c_2^{H,\gamma} =c_2^{H}$. Furthermore, for any abelian group there is a
 canonical  homomorphism
   \[ l_{2}(A)\:\colon\hspace{2mm} A \sm A \to A \ot A\:,\quad l_{2}(A)(x\sm y) = x\ot y - y\ot x\]
for $x,y\in A$. We note $l_{2}(G/G_{(2)}) = l_{2}^{\cal G}$,   $l_{2}(H/H_2) = l_{2}^{H}$ and $l_2^{\calG H} = (\iota \ot id)l_2^H$
with $\iota\,\colon\, H/H_2 \to G/KG_{(2)}$ induced by $H\hra G$.
Now consider the
commutative square
 \BE\label{squFox2} \begin{matrix}
\ruleu\ruled H\ab \sm H\ab   &  \mr{l_2^{\calG H}}  &  (G/KG_{(2)}) \ot (H/H_2) \cr
\surdown{c_2^{H}}  &  & \surdown{\mu_2^{\calG H}} \cr
\rule{0mm}{6mm} H_2/H_3  &  \mr{d_2^{\calG H}}  &    I(G) I(H) \Big/  \Big(I(K)I(H) + I^2_{\calG}(G) I(H)
\Big)\rule[-5mm]{0mm}{3mm}  \rule[-5mm]{0mm}{3mm}\end{matrix}\EE
\rule{0mm}{7mm} where for $g \in G$, $h \in H$  and $x \in H_2$
  \[ \mu_2^{\calG H}((gG_{(2)})\ot (hH_2)) = \overline{(g-1)(h-1)} \:, \quad  d_2^{\calG H}(xH_3) = \overline{x-1}\:. \]

This diagram provides a simultaneous functorial description of both the second Fox
quotient and the second Fox subgroup, in a completely symmetric way, as
follows.\vspace{0.3mm}

\begin{satz}\label{Fox2IH}\quad Diagram \REF{squFox2} is a pushout square of abelian
groups; in particular
  \[ I(G)I(H)\Big/\Big(I(K)I(H) + I^2_{\calG}(G) I(H)\Big) \:\cong\: (G/KG_{(2)}) \ot (H/H_2) \Big/ l_2^{\calG H}  {\rm
Ker}(c_2^{H}) \]
\[  \Big(G \cap \Big(1 + I(K)I(H) + I^2_{\calG}(G)I(H)\Big) \Big) \Big/  H_3 = \Ker{d_2^{\calG H}} = c_2^{H}\,\Ker{l_2^{\calG
H}}\:.\]
\end{satz}\vspace{0mm}

An analogous result holds for the third Fox quotient and subgroup, see Theorem
\ref{Fox3push} below.\V

\N{\bf Proof\,:}\quad In Proposition \ref{satz} take $n=2$, $J=I_R(H)$ and $U=I^2(H)$. Then using the isomorphisms $\phi_H$ and
$\phi_G$ we see that the following diagram is a pushout square, with ${\cal H}=\gamma_H$.
  \BE\label{Fox2pushR} \begin{matrix}
\ruleu\ruled H\ab \ot H\ab \ot R & \mr{\iota \ot id  \ot id}  & ( G/KG_{(2)}) \ot H\ab  \ot R \cr
\surdown{\mu_2^{{\cal H} H} \ot id }  &  &  \surdown{\mu_2^{\calG H} \ot id } \cr
I_R^2(H)/I_R^3(H)  &  \mr{j}  &  I_R(G)I_R(H) \Big/ \Big(I_R(K)I_R(H) + I_{R,\cal G}^2(G) I_R(H) \Big) \ruleu\rule[-5mm]{0mm}{3mm}
\end{matrix}\EE
So $\Ker{\mu_2^{\calG H}\ot id} = (\iota \ot id\ot id) \Ker{\mu_2^{{\cal H} H}\ot id}$. But for $R=\Z$, $\Ker{\mu_2^{{\cal H} H}} 
= l_2^{H} \Ker{c_2^H}$ by the isomorphism
$I^2(H)/I^3(H) \:\cong\: {\rm U}_2{\rm
L}(H)$ obtained in \cite{Ba-GrII}. Hence diagram \REF{squFox2} is a pushout, too (abstractly, this follows from the gluing lemma
for pushouts in any category as the cited isomorphism means that  diagram \REF{squFox2} for $K=\{1\}$, $G=H$ and ${\cal G}=\gamma$ is a pushout
square). 
To deduce from this fact the second identity in the assertion, we first use Corollary \ref{GcapIJprop} to observe that
$G \cap (1 + I (G)I(H)) = H \cap (1 + I(H)^2)=H_2$.\hfbox\V

We remark that in \cite{D3F2} we deduce a more explicit description of the induced subgroup $G\cap (1  +
I_R(K)I_R(H) + I_R^2(G) I_R(H))$ from the pushout square \REF{Fox2pushR}, for
arbitrary subgroups $K,H$ of $G$ and coefficient rings $R$.

A nice generalization of Corollary \ref{Fox2IH} also holds  for higher Fox quotients in the
case where $H/H_n$ is free nilpotent which morally means that $H$ does not have any relations in lower commutator filtration. This is
the matter of the next section.\V

\section{Fox quotients with respect to free nilpotent subgroups} 

In Proposition \ref{Hfreequot} we saw that the map $\zeta_n^{\calG}$ is an isomorphism if $H$ is free.
In this section we compute its kernel for $\calG =\gamma$ if only $H/H_n$ is  free nilpotent of class $n-1$.\V

Let ${\cal L}(H\ab)$ and $T(H\ab) = \Z \oplus \bigoplus_{n>0} (H\ab)\htt{n}$  be the
free Lie ring and the tensor ring over $H\ab$, respectively. Let L$(H)$ be the graded Lie
ring defined by the successive lower central quotients of $H$ and \,Gr$(\Z(H))$ be the
associated graded ring of $\Z(H)$ with respect to the filtration $I^k(H)$, $k\ge 0$. Then we
have natural  homomorphisms of graded abelian groups (actually, of graded rings in the case of $q$ and of graded Lie rings in the
case of $l$ and $c$)
  \[ {\rm Gr}(\Z(H)) \surltop{q} T(H\ab) \ml{l} {\cal L}(H\ab) \Sur{c} {\rm
L}(H)  \] 
provided by the identifications $q_1 = l_1 =  c_1 = id_{H\ab}$ and the universal
properties of $T(H\ab)$ and ${\cal L}(H\ab)$. If $H$ is free then $c$ is well-known to
be an isomorphism. So for every $n\ge 2$ we obtain natural homomorphisms
  \[  Q_{n-1}(G)\ot H\ab \;\ml{Q_{n-1}(i)q_{n-1} \ot id}\; (H\ab)\htt{n} \ml{l_{n}}
{\cal L}_{n}(H\ab) \Sur{c_{n}^H}  {\rm L}_{n}(H)  \]

 The maps $l_n$ and $c_n$ may be viewed as sending a ``formal''
$n$-fold commutator in ${\cal L}_{n}(H\ab)$ to the corresponding tensor commutator and
to the coset of the corresponding group commutator, respectively. With these notations, we
get the following

\begin{satz}\label{Foxfreequot}\quad Let $n\ge 2$ and suppose that $H/H_n$ is a free
nilpotent group of class $n-1$, i.e.\ $H/H_n \cong F/F_n$ for some free group $F$. Then the
map $\zeta_n$  above induces a natural isomorphism
  \[ I^{n-1}(G)I(H)/I^{n}(G)I(H) \:\cong\: Q_{n-1}(G)\ot H\ab/ (Q_{n-1}(i)q_{n-1} \ot
id) l_n\Ker{c_n^H}\:.\]
In particular, if  $H/H_{n+1}$ is free nilpotent of class $n$, then 
  \[ I^{n-1}(G)I(H)/I^{n}(G)I(H) \:\cong\: Q_{n-1}(G)\ot H\ab \:.\]
\end{satz}

\N We remark that the proof below shows that the group
$(Q_{n-1}(i)q_{n-1} \ot id) l_n\Ker{c_n}$ is contained in the kernel of $\zeta_n$
for {\em arbitary}\/ subgroups $H$; it actually equals ${\rm Ker}(\zeta_n)$ for $n=2$, by Theorem \ref{Fox2IH}. In
the remaining sections we study a refinement of the approximation of the $n$-th Fox quotient by $\zeta_n$ which takes  this fact into account (among other phenomena).\V

\proof We wish to apply Proposition \ref{satz} with $U=I^n(H)$. By  assumption on $H$ there is a free
presentation $\,R \hra F \Sur{\pi} H$ such that $R \subset F_n$. It induces a ring
isomorphism $\Z(F)/I^k(F) \,\cong\,\Z(H)/I^k(H)$ for all $k\le n$. Choose a basis
$(f_i)_{i\in I}$ of $F$. Since
$I(F)$ is a free $F$-module with basis $(f_i-1)_{i\in I}$, it follows that $I(H)/I^n(H)$ is a 
free $\Z(H)/I^{n-1}(H)$-module with basis $(\pi(f_i)-1)_{i\in I}$. Whence in diagram
\REF{diasatz} (with $J=I(H)$) we obtain an isomorphism
  \BE\label{xi} \xi\,\colon\,  P_{n-1}(G) \ot H\ab  \:\stackrel{\cong}{\to}\:
P_{n-1}(G) \ot_H P_{n-1}(H) \EE 
defined by $\xi( x\ot \pi(f_i)H\st) = x \ot (\pi(f_i)-1+I^n(H))$.
Since $H\ab$ is free abelian we get an isomorphism
  \BE\label{Habfrei} \xi(\iota\ot id)\,\colon\,Q_{n-1}(G)\ot H\ab
\:\cong\: \Imm{( Q_{n-1}(G)  \ot_H P_{n-1}(H) }
  \EE
with the inclusion $\iota\,\colon\, Q_{n-1}(G) \hra  P_{n-1}(G)$.
 In order to compute $\Ker{\mu_H}$ in \REF{diasatz} write $T= T(H\ab)$ and let $\bar{T}$ be its
augmentation ideal. 
 It is well known that there is a 
ring isomorphism $\chi\,\colon\, T/\bar{T}^k \:\cong \: \Z(F)/I^k(F)$ for all $k\ge 1$,
defined by sending a generator $f_iF_2$ to $f_i-1+ I^k(F)$. Thus we get a commutative
square
  \[\begin{matrix}
\ruled P_{n-1}(H)  \ot_H  P_{n-1}(H)  & \isol{ \chi\ot \chi} &
(\bar{T}/\bar{T}^n) \ot_T  (\bar{T}/\bar{T}^n)  & \isor{\nu} & \bar{T}^2/\bar{T}^{n+1}  \cr
\mapdown{\mu_H} & & & & \surdown{}\cr   
\ruleu\frac{\dst I^2(H)}{\dst I^{n+1}(H)} & \cong & \frac{\dst
I^2(F)+I(R)\Z(F)}{\dst  I^{n+1}(F)+I(R)\Z(F)}  &  \cong & \bar{T}^2/(\bar{T}^{n+1} +
p_n\Ker{\pi_n})
\end{matrix}\]
where $\nu$ is induced by multiplication in $T$ and where the maps $p_n,\pi_n$ are part of the following commutative diagram and are
induced by sending $a \in F_n$ to $\chi^{-1}(a-1 +I^{n+1}(F))$ and to $\pi(a)H_{n+1}$,
respectively. The isomorphisms in the bottom row above are then deduced from the relations
$I(R) \subset I(F_n) \subset I^n(F)$. \\
\N\makebox[14.7cm]{ \makebox[0mm]{
\begin{minipage}{20cm}\small
\[ \begin{matrix}
\ruleu (\bar{T}/\bar{T}^n) \ot_T  (\bar{T}/\bar{T}^n) & \isol{\nu} & \bar{T}^2/\bar{T}^{n+1}
&\ml{p_n}  &  F_n/F_{n+1}  & \Sur{\pi_n}  & H_n/H_{n+1}  \cr 
\isodown{\chi\ot\chi}  &  & \mapup{\omega}  &  &
\cong\hspace{-1mm}\mapup{c_n^F} 
 &  &  \mapup{c_n^H} \cr
\ruled P_{n-1}(H) \ot_H  P_{n-1}(H)   &   &
(H\ab)\htt{n}  &  \ml{l_n{\cal L}_n(\pi\ab)}  &  {\cal L}_n(F\ab) & \isor{{\cal
L}_n(\pi\ab)}   &  {\cal L}_n(H\ab)  \cr
\mapdown{i\ot id} & & \mapdown{(Q_{n-1}(i)q_{n-1}
\ot id)} & & & & \cr
\rule{0mm}{5mm}  P_{n-1}(G)  \ot_H  P_{n-1}(H)   & \ml{\xi( \iota\ot id)}  & 
Q_{n-1}(G)\ot H\ab & & & &   \end{matrix}\]
\end{minipage}\ruled
}  }
\rule{0mm}{7mm}

\N Here $\omega$ is the canonical injection. From the first of these two diagrams it follows that $\Ker{\mu_H} =
(\chi\ot\chi)\nu^{-1} p_n
\Ker{\pi_n}$, so by \ref{satz} and the second diagram 
$\Ker{\mu_G} = (i\ot id)(\chi\ot\chi)\nu^{-1}p_n\Ker{\pi_n} = (i\ot id)(\chi\ot\chi)\nu^{-1}\omega l_n\Ker{c_n^H} = 
\xi(\iota\ot id)
(Q_{n-1}(i)q_{n-1} \ot id) l_n\Ker{c_n^H}$. By the equality $\mu_G \xi( \iota\ot H\ab) =
\zeta_n$ and by \REF{Habfrei} this implies the
assertion. \hfbox

\section{Canonical approximation of Fox quotients}


In the preceding  sections we computed the generalized Fox quotients $ Q_{n}^{\cal G}(G,H)$ in some special cases by using the
somewhat ``naive" approximation $\zeta_n^{\calG}$; we here study the  deeper structure of the groups $ Q_{n}^{\cal G}(G,H)$ by
introducing a much closer approximation in terms of enveloping algebras which generalizes
Quillen's
  approximation of augmentation quotients, i.e. the case $\calG =\gamma$ and $H=G$. We start by recalling the latter construction.

The abelian group ${\rm L}^{\cal G}(G) = \sum_{n\ge 1} G_{(n)}/G_{(n+1)}$ is a graded Lie ring whose bracket is induced by the
commutator pairing of
$G$. So its enveloping algebra $\ULG{}$ over the integers is defined. On the other hand, the  filtration quotients
$Q_n^{\calG}(G) = I^n_{\cal G}(G)/I^{n+1}_{\cal G}(G)$ form the graded ring  Gr$^{\cal G}(\Z(G)) = \oplus_{n\ge 0} Q_n^{\calG}(G)$;
note that one has  Gr$(\Z(G)):\,=\,$Gr$^{\gamma}(\Z(G))=  \bigoplus_{n\ge 0} I^n(G)/I^{n+1}(G)$. 


Now the map 
$\LG{} \to$ Gr$^{\cal G}(\Z(G))$, $aG_{(n+1)} \mapsto a-1+I^{n+1}_{\calG}(G)$ for $a\in G_{(n)}$, is a homomorphism of graded Lie
rings and hence extends to a map of graded rings $\theta^{\calG}\,\colon\,\ULG{} \to$ Gr$^{\cal G}(\Z(G))$. This map is clearly
surjective but rarely globally injective; for instance, $\theta^{\gamma}$ is  
injective  if $G$ has torsionfree lower central quotients $G_n/G_{n +1}$ or is cyclic, but $\theta^{\gamma}$ is non
injective for all non cyclic finite abelian groups \cite{Ba-GrI}. At least, the kernel of  $\theta^{\calG}$ is torsion as
$\theta^{\calG} \ot \Q$ is an isomorphism; this was proved by
 Quillen for ${\cal G}=\gamma$ and  follows from  work of Hartley \cite{Hartley} in the general case, see
also \cite{Foxsplit}. Moreover, \Ker{\theta^{\calG}} is trivial in degree 1 and 2 (by \cite{Ba-GrII} for $N=\gamma$) and is
explicitely known in degree 3, see \cite{Q3}.


To generalize the foregoing concept to Fox quotients consider the  filtration 
\[ {\cal F}^1 = \Z(G)I(H) \supset {\cal F}^2 =I(G)I(H) \supset \ldots \supset
{\cal F}^n =I^{n-1}_{\cal G}(G)I(H) \supset \ldots \]
 of $\Z(G)I(H)$. 
%
%
 The associated graded group \,Gr$^{{\cal G}\gamma}(\Z(G)I(H)) =
\bigoplus_{n\ge 1} Q_n^{\cal G}(G,H)$ is  a  graded Gr$^{\cal G}(\Z(G))$--Gr$(\Z(H))$--bimodule in the canonical way, and hence a
$\ULG{}-{\rm UL}(H)$-bimodule via the maps
$\theta^{\calG}$ and $\theta^{\gamma_H}$. Now let ${\cal H}= (H_{(n)})_{n\ge 1}$  be the N-series of
$H$ defined by $ H_{(n)} = H\cap G_{(n)}$.
 The injection $I(H)
\hra \Z(G)I(H)$ is a   map of $\Z(H)$-bimodules taking $I^n(H)$ and $I^{n-1}_{\cal H}(H)I(H)$ into ${\cal F}^n$;
  it thus induces   homomorphisms
\,Gr$(I(H))
\mr{}$ Gr$^{{\cal G}\gamma}(\Z(G)I(H)) \:\ml{}\:$ Gr$^{{\cal H}\gamma}(I(H))$ of graded \,Gr$(\Z(H))$-bimodules and \,Gr$^{\cal
H}(\Z(H))$--Gr$(\Z(H))$--bi\-mo\-dules, resp. So by extension of scalars along the graded ring homomorphisms ${\rm Gr}(\Z(H))
\to {\rm Gr}^{\cal G}(\Z(G))$ and ${\rm Gr}^{\cal H}(\Z(H))
\to {\rm Gr}^{\cal G}(\Z(G))$ 
induced by the injection $H\hra G$ we get  natural surjective maps of graded
${\rm Gr}^{\cal G}(\Z(G))$--Gr$(\Z(H))$--bimodules
   \[   \xi^{GH}_{{\cal G}\gamma}\:\colon\hspace{2mm}  {\rm Gr}^{\cal G}(\Z(G)) \ot_{\raisebox{-1mm}{$\sst {\rm Gr}(\Z(H))$}} {\rm
Gr}(I(H)) 
\:\mr{}\: {\rm Gr}^{\calG \gamma}(\Z(G)I(H)) \:,\]
  \[   \xi^{GH}_{ \cal GH} \:\colon\hspace{2mm}  
{\rm Gr}^{\cal G}(\Z(G)) \ot_{\raisebox{-1mm}{$\sst {\rm Gr}^{\cal H}(\Z(H))$}} {\rm
Gr}^{{\cal H}\gamma}(I(H)) \:\mr{}\: {\rm Gr}^{{\cal G}\gamma}(\Z(G)I(H)) \:.\]
It is convenient to
combine $\xi^{GH}_{{\cal G}\gamma}$ with Quillen's approximation to obtain an epimorphism of graded $\ULG{}-{\rm UL}(H)$-bimodules
  \[ \theta^{\calG H} = \xi^{GH}_{{\cal G}\gamma} (\theta^{\calG} \ot \theta^{\gamma_H}) \hspace{1mm}\colon \hspace{2mm} \UGH{} \hspace{1mm}:\, =
\hspace{1mm} \ULG{}
\ot_{\raisebox{-1mm}{$\sst {\rm UL}(H)$}}
{\rm \bar{U}L}(H)
\:\mr{}\: {\rm Gr}^{{\cal G}\gamma}(\Z(G)I(H)) \]
where ${\rm \bar{U}L}(H)$ denotes the augmentation ideal of ${\rm UL}(H)$.

\begin{bem}\rm The approximation of ${\rm Gr}^{{\cal G}\gamma}(\Z(G)I(H))$ by $\xi^{GH}_{ \cal GH}$ is a priori ``closer'' than the
one by 
$\xi^{GH}_{{\cal G}\gamma}$ as is seen from the following commutative diagram
\\
\N\makebox[14.7cm]{ \makebox[0mm]{
\begin{minipage}{20cm}\small
  \[ \begin{matrix}
\ruleu\rule[-5mm]{0mm}{0.5mm} {\rm Gr}^{\calG}(\Z(G))  \ot_{\raisebox{-1mm}{$\sst {\rm Gr}(\Z(H))$}} {\rm Gr}(I(H))  & \isor{m}  & 
{\rm Gr}^{\calG}(\Z(G))   \ot_{\raisebox{-1mm}{$\sst {\rm Gr}^{\cal H}(\Z(H)) $}}  \Big( {\rm Gr}^{{\cal H} }(\Z(H)) 
\ot_{\raisebox{-1mm}{$\sst {\rm Gr}(\Z(H))$}} {\rm Gr}(I(H))  \Big)  \cr
\mapdown{\xi^{GH}_{ {\cal G}\gamma } } &  &  \mapdown{id\ot \xi^{HH}_{{\cal H}\gamma}} \cr
{\rm Gr}^{{\cal G}\gamma}(\Z(G)I(H))  &  \ml{\xi^{GH}_{{\cal GH}} }  &  
{\rm Gr}^{\calG}(\Z(G))  \ot_{\raisebox{-1mm}{$\sst {\rm Gr}^{\cal H}(\Z(H))$}} 
\rule{0mm}{5mm}\ruled {\rm Gr}^{{\cal H}\gamma}(I(H)) 
\end{matrix} \rule[-5mm]{0mm}{3mm} \]
\end{minipage}\ruled
}}
 \vspace{4mm}

\N where $m$ is the canonical isomorphism. But as our goal is to approximate the group ${\rm Gr}^{{\cal G}\gamma}(\Z(G)I(H))$ in terms of 
enveloping algebras we do not to care about this difference.
\end{bem}

 For $i\ge 0$ and $j\ge 1$ let 
   \[  \nu_{ij}\,\colon\,\ULG{i} \ot {\rm U}_j{\rm L}(H) \to
{\rm U}_{i+j}^{{\cal G}}(G,H)  \]
 be the canonical map. 
Note that for any group $K$ the ring ${\rm U}{\rm L}(K)$ is generated by
${\rm U}_1{\rm L}(K) = {\rm L}_1(K) \:\cong\:K\ab$. This implies that $\nu_{(n-1)1}$ is surjective, and also implies  exactness of the
following sequence of graded
$\ULG{}-{\rm UL}(H)$-bimodules
  \[ \ULG{} \ot {\rm U}_1{\rm L}(H) \ot {\rm \bar{U}L}(H) \hspace{2mm}\mr{\psi}\hspace{2mm} \ULG{} \ot {\rm \bar{U}L}(H) \hspace{2mm}\mr{q}\hspace{2mm}
\UGH{} \hspace{2mm}\to 0 \] where $\psi(x\ot y \ot z) = xy \ot z - x \ot yz$. It is now easy to compute $\UGH{n}$ for $n\le 3$,
also using  the identity
$\Ker{\mu_2^{H}} = l_2^{H}\,\Ker{c_2^{H}}$, cf.\ the proof of Theorem \ref{Fox2IH}.

From now on we abbreviate $G^{AB} =G/G_{(2)}$.\V

\begin{prop}\label{UnGH} There are canonical isomorphisms
  \[ \UGH{1} \hspace{2mm}\cong\hspace{2mm}   {\rm U}_1{\rm L}(H)  \hspace{2mm}\cong\hspace{2mm}  H\ab  \]
  \[ \UGH{2} \hspace{2mm}\cong\hspace{2mm}   \frac{\dst G^{AB} \ot H\ab}{\ruled \dst l_2^{\calG H}\Ker{c_2^{H}}}  \]
\[ \UGH{3} \hspace{2mm}\cong\hspace{2mm} {\rm coker}\Big( \epsilon = \left( 
\begin{array}{rcc}
 c_2^{\calG} \ot id &   0  &    0  \cr
  {}-l_2^{\calG} \ot id  &   i^{\calG HH}   &    i^{HHH} 
\end{array}\right) \Big)\:, \]

  \[\begin{matrix}
\ruled
( (G^{AB} \sm G^{AB}) \ot H\ab ) \hspace{4mm}\oplus\hspace{4mm} ( G^{AB} \ot  l_2^H \Ker{c_2^{H}} )
\hspace{4mm}\oplus\hspace{4mm}
l_{33}^H \Ker{c_{33}^H} \cr
\mapdown{\epsilon }  \cr
\rule{0mm}{6mm}\ruled
(  ( G_{(2)}/ G_{(3)} ) \ot H\ab )  \hspace{4mm}\oplus\hspace{4mm}  ( G^{AB} \ot G^{AB} \ot H\ab )
\end{matrix}\]

\N where the homomorphisms ${\rm L}_3(H)  \ml{c_{33}^H}  (H\ab)^{\ot 3}  \mr{l_{33}^H}  (H\ab)^{\ot 3}$ are defined such that for
$x,y,z\in H\ab$, 
$c_{33}^H(x\ot y
\ot z)$ is the triple Lie bracket $[x,[y,z]]$ in the Lie algebra ${\rm L}(H)$ and  $l_{33}^H(x\ot y \ot
z)$ is the triple Lie bracket $[x,[y,z]]$ in the tensor algebra $T(H\ab)$. Furthermore, we note $i^{\calG HH} =  id \ot \iota \ot
id\,\colon\,G^{AB} \ot H\ab \ot H\ab \to G^{AB} \ot G^{AB} \ot H\ab $ 
 and
$i^{HHH} = \iota \ot \iota \ot id\,\colon\, H\ab \ot H\ab \ot H\ab  \to  G^{AB} \ot G^{AB} \ot H\ab $.
\end{prop}\vspace{2mm}

As an immediate consequence of this computation and of Proposition \ref{Fox1} and Theorem \ref{Fox2IH} we get

\begin{prop}\label{xi1et2} For all groups $G$ and subgroups $H$ of $G$ the maps $\xi_n^{\calG H}$ and $\theta_n^{\calG H}$ are
isomorphisms for $n=1,2$.
\end{prop}


One may also ask when  $\theta^{\calG H}$ is globally an isomorphism. We have a positive answer in at least one case.\V

\begin{kor}\label{thetaHfree} If $\calG$ has torsionfree factors and $H$ is a free group then $\theta^{\calG H}$ is an 
isomorphism.
\end{kor}


\N{\bf Proof\,:} We have the following commutative diagram
  \[\begin{matrix}
\ruleu\ruled {\rm U}_{n-1}\LG{} \ot H\ab & \mr{\theta_{n-1}^{\calG} \ot id}  & Q^{\calG}_{n-1}(G) \ot H\ab \cr
\surdown{\nu_{(n-1)1}}  &  &  \isodown{\zeta_n^{\calG}} \cr
\ruleu\ruled  \UGH{n}  &  \Sur{\theta_n^{\calG H}}  &  {\rm Gr}_n^{\calG \gamma}(\Z(G)I(H)) 
\end{matrix}\]
where $\zeta_n^{\calG}$ is an isomorphism by Proposition \ref{Hfreequot}.
Moreover, $\theta^{\calG}$ here is an isomorphism since  $\calG$ has torsionfree factors \cite{Foxsplit}, whence $\nu_{(n-1)1}$ and 
$\theta_n^{\calG H}$ are isomorphisms, too.\hfbox\V

We now exhibit a canonical part of the kernel of $\theta^{\calG H}$ (other than \Ker{\theta^{\calG} \ot \theta^{\gamma_H}}) 
in comparing $\calG$ to the lower central series of
$H$, as follows.

For  elements $x_1,\ldots,x_m$ of any ring define the  iterated  commutator  $[x_1,\ldots,x_m]$ to be  $x_1$
if
$m=1$ and 
to be  $[x_1,[x_2,\ldots,[x_{m-1},x_m]\ldots]$ if $m\ge 2$. In the latter case one has the formula
   \BE\label{commformel} [x_1,\ldots, x_m] =
\sum_{s=0}^{m-1}(-1)^{m-s+1} \sum_J x_{i_1} \cdots x_{i_s} x_m x_{j_1}\cdots
x_{j_{m-1-s}}  \EE
 where $J= \{(i_1,\ldots,i_s,j_1,\ldots, j_{m-1-s})\,\|\,  1\le i_1
<\ldots< i_s \le m-1  \ge j_1  >\ldots>j_{m-1-s} \ge 1 \mbox{ such that}
 \{i_1,\ldots,i_s,j_1,\ldots,  j_{m-1-s}\} = \{1,\ldots,m-1\} \,\} $. 
We say that an $m$-tuple
 ${\underline{h}} = (h_1,\ldots,h_m)$ of elements of $H$ is of {\em height}\/ 
 $\ge n$ if $h_j \in H_{k_j} \cap G_{({l_j})}$ such that $l_1 +\ldots + l_m - l_j + k_j = n$ for $1 \le j \le m$. 
For such an
$m$-tuple ${\underline{h}}$ let $r_1 ({\underline{h}}) = h_1$ if $m=1$ and if  $m\ge 2$,
  \[ r_1 ({\underline{h}}) =  [h_1,[h_2, \ldots,[h_{m-1},h_m]\ldots]   \in H \]
  \[  r_2 ({\underline{h}}) =  \sum_{s=0}^{m-1}(-1)^{m-s+1} \sum_J  (h_{i_1} G_{({l_{i_1} +1})}) \ldots (h_{i_s} G_{({l_{i_s} +1})}). 
(h_{m} G_{({l_{m} +1})}). (h_{j_1} G_{({l_{j_1} +1})}) \ldots 
\]
  \[ \Ph\hspace{21mm} (h_{j_{m-2-s}} G_{({l_{j_{m-2-s}} +1})})   \ot  (h_{j_{m-1-s}} H_{k_{j_{m-1-s}} +1}) 
\in  \Big( \ULG{}\ot_{{\rm UL}(H)} {\rm \bar{U}L}(H) \Big)_n \:\]
\vspace{1mm}

\N Here, and throughout the rest of this paper,  we consider the cosets $h_jH_{k_j+1} \in {\rm L}_{k_j}(H)$ and $h_jG_{({l_j +1})} \in {\rm L}^{\calG}_{l_j}(G)$ as elements
of
$ {\rm U}_{k_j}{\rm L}(H)$ and
$ {\rm U}_{l_j}\LG{}$, resp., suppressing the canonical map ${\rm L}^{\cal K}(K) \to {\rm U}{\rm L}^{\cal K}(K)$ from the notation. Moreover,  all products denoted by $.$ are taken in the respective enveloping algebras.

For $n\ge 2$ let ${\cal R}_n^{\calG H}$  
be the subgroup of $\ULG{}\ot_{{\rm UL}(H)} {\rm
\bar{U}L}(H) $   generated by the elements 
  \[ R_n({\underline{h}}_1,\ldots,{\underline{h}}_p) = 1_{\ULG{}} \ot (r_1({\underline{h}}_1) \cdots
r_1({\underline{h}}_p)   H_{n+1}) - \sum_{q=1}^p r_2({\underline{h}}_q)   \in   {\rm U}_n^{\calG}(G,H)
\] 
where $p\ge 1$ and each ${\underline{h}}_q$ is an $m_q$-tuple, $m_q\ge 2$, of height $\ge n$ such that 
 $r_1({\underline{h}}_1) \cdots r_1({\underline{h}}_p) \in H_n$.   The term ${\cal R}_3^{\calG H}$ will be considered in Corollary \ref{R3=} below.\V

\begin{prop}\label{kerneltheta}\quad One has identities $\theta^{\calG H}({\cal R}_n^{\calG H}) = 0$ and ${\cal R}_n^{\calG H} \, {\rm \bar{U}L}(H) =0$.
Moreover, in the definition of ${\cal R}_n^{\calG H}$ it suffices to take only those tupels ${\underline{h}}_q = (h_{q1},\ldots,h_{qm_q})$ for which $l_{qj} \ge
k_{qj}$ for $1 \le j \le m_q$. 
\end{prop}\V

\proof Let ${\underline{h}} = (h_1,\ldots,h_m) \in H^m$ be of height $\ge n$. We contend that in $\Z(G)$, $r_1({\underline{h}}) - 1 \in {\cal
F}^n $ and if $m\ge 2$, 
  \BE\label{r1(H)}   
  r_1({\underline{h}}) - 1
\:\equiv\: [h_1-1,[\ldots,[h_{m-1} -1, h_m -1]\ldots ] \hspace{1mm}\mbox{mod}\hspace{1mm} {\cal F}^{n+1}\:. \EE

\N We proceed by induction on $m$. The case $m=1$ being obvious let $m\ge 2$. Note that ${\underline{h}}\st = (h_2,\ldots,h_m) $ is of
height $\ge n-l_1$; hence by induction hypothesis, $r_1({\underline{h}}\st) - 1 \in {\cal F}^{n-l_1}$ and $x= r_1({\underline{h}}\st) - 1 - 
[h_2-1,[\ldots,[h_{m-1} -1, h_m -1]\ldots ] \in {\cal F}^{n-l_1+1}$. On the other hand,   $r_1({\underline{h}}\st) - 1 \in 
\theta^{\calG} (\calG_{l_2+\ldots+l_m}) \subset
I_{\calG}^{l_2+\ldots+l_m}(G)$ and that $x\in I_{\calG}^{l_2+\ldots+l_m+1}(G)$ since $\theta^{\calG}$ is a graded Lie map.
Now
 \[  r_1({\underline{h}}) - 1  =   \Big( (h_1 -1)(r_1({\underline{h}}\st) - 1)  -  (r_1({\underline{h}}\st) - 1)(h_1 -1) \Big)h_1^{-1} r_1({\underline{h}}\st)^{-1} \:. \]
By the above remarks, $(h_1 -1)(r_1({\underline{h}}\st) - 1)  -  (r_1({\underline{h}}\st) - 1)(h_1 -1) \in I_{\calG}^{l_1}(G){\cal F}^{n-l_1}  + 
I_{\calG}^{l_2+\ldots+l_m}(G) I^{k_1}(H) \subset {\cal F}^n$ by the immediate relation $I_{\calG}^l(G) {\cal F}^e I^k(H) \subset
{\cal F}^{l+e+k}$ for $l,e,k\ge 0$. Thus also $ r_1({\underline{h}}) - 1\in {\cal F}^n$. Writing $h_1^{-1} r_1({\underline{h}}\st)^{-1} = 1 + (h_1^{-1} r_1({\underline{h}}\st)^{-1} -1)$ we get 
\begin{eqnarray*}
r_1({\underline{h}}) - 1  & \equiv  &  (h_1 - 1)(r_1({\underline{h}}\st) - 1) -  (r_1({\underline{h}}\st) - 1)(h_1 - 1)\quad \mbox{mod ${\cal F}^{n+1}$}\\
  &  =  & [h_1 - 1\,, x + [h_2-1,[\ldots,[h_{m-1} -1, h_m -1]\ldots ]  ] \\
  &  \equiv  & [h_1 - 1\,,[h_2-1,[\ldots,[h_{m-1} -1, h_m -1]\ldots ]] \quad \mbox{mod ${\cal F}^{n+1}$}
\end{eqnarray*}
since $(h_1 - 1)x - x(h_1 - 1) \in I_{\calG}^{l_1}(G) {\cal F}^{n-l_1+1}  +  I_{\calG}^{l_2+\ldots+l_m+1}(G) I^{k_1}(H)  \subset
{\cal F}^{n+1}$. Hence \REF{r1(H)} is proved. Using the identity $ab - 1  =  (a-1)+(b-1)+(a-1)(b-1)$ we obtain
\begin{eqnarray*}
\theta^{\calG H}(1_{\ULG{}} \ot r_1({\underline{h}}_1) \ldots r_1({\underline{h}}_p) H_{n+1})  
  &=&  r_1({\underline{h}}_1) \ldots r_1({\underline{h}}_p) -1 + {\cal F}^{n+1} \\
  &=&  \sum_{q=1}^p r_1({\underline{h}}_q) -1 + {\cal F}^{n+1} \\
  &=&  \sum_{q=1}^p \theta^{\calG H} (r_2({\underline{h}}_q)  )+ {\cal F}^{n+1}  
\end{eqnarray*}
where the last identity is due to relations \REF{r1(H)} and  \REF{commformel}. Hence $\theta^{\calG H}({\cal R}_n^{\calG H}) = 0$. Next we prove the last part
of the assertion. First note that for $\underline{h}$ as above
\BE\label{l1+lm}
l_1+\ldots+l_m = n+l_j-k_j
\EE
for $1\le j\le m$. Now  suppose that for some $q$, $l_{q1}+\ldots+l_{qm_q} < n$ which by \REF{l1+lm} means that  $l_{qj} < k_{qj}$ for {all}\/ $j$
from
$1$ to
$m_q$. For clarity we suppress the index $q$ from  the notation. We have $h_j \in H_{k_j} \subset G_{({k_j})} \subset G_{({l_j+1})}$, so
$r_2({\underline{h}}_q) = 0$ since all terms to the left of $\ot$ are trivial. On the other hand, 
  \begin{eqnarray*}
k_1+ \ldots + k_m  &=& n - (l_2 + \ldots +l_m) + k_2 +\ldots + k_m \\
  &=&  n + (k_2-l_2)  +\ldots + (k_m-l_m) \\
  &\ge& n + m-1\\
  &\ge& n + 1 \:,
\end{eqnarray*}
so $r_1({\underline{h}}_q) \in H_{k_1   +\ldots + k_m} \subset H_{n+1}$. Hence  $R_n({\underline{h}}_1,\ldots,{\underline{h}}_p) =  R_n({\underline{h}}_1,\ldots,\hat{\underline{h}}_q,\ldots, {\underline{h}}_p)$ where $({\underline{h}}_1,\ldots,\hat{\underline{h}}_q,\ldots, {\underline{h}}_p)$ is $({\underline{h}}_1,\ldots,{\underline{h}}_p)$ with ${\underline{h}}_q$ omitted. Thus we may suppose
that in the definition of ${\cal R}_n^{\calG H}$, for all $q$ from $1$ to $p$, one has
$l_{q1} +\ldots + l_{qm_q} \ge n$, i.e., $l_{qj} \ge
k_{qj}$ for $1 \le j \le m_q$ by \REF{l1+lm}. Under this hypothesis let $v\in  \bar{\rm U}{\rm L}(H)$ and let us  show that $R_n({\underline{h}}_1,\ldots,{\underline{h}}_p)
\cdot v=0$. We have
  \[  R_n({\underline{h}}_1,\ldots,{\underline{h}}_p) \cdot v   =
 (r_1({\underline{h}}_1) \cdots
r_1({\underline{h}}_p)   G_{({n+1})}) \ot v  \hspace{2mm} - \hspace{2mm} \sum_{q=1}^p r_2^G({\underline{h}}_q)  \ot v 
\] 
where 
  \[  r_2^G ({\underline{h}}) =  \sum_{s=0}^{m-1}(-1)^{m-s+1} \sum_J  (h_{i_1} G_{({l_{i_1} +1})}) \ldots (h_{i_s} G_{({l_{i_s} +1})}). 
(h_{m} G_{({l_{m} +1})}) .(h_{j_1} G_{({l_{j_1} +1})}) \ldots 
\]
  \[ \Ph\hspace{26mm} (h_{j_{m-2-s}} G_{({l_{j_{m-2-s}} +1})})  .   (h_{j_{m-1-s}} G_{({k_{j_{m-1-s}} +1})}) 
\]\vspace{1mm}
 
\N Suppose that for some $q$, $l_{q1} +\ldots + l_{qm_q} > n$, i.e., $l_{qj} >
k_{qj}$ for $1 \le j \le m_q$. Then $h_{qj} G_{({k_{qj}+1})} = 0$ as  
$h_{qj} \in G_{({l_{qj}})} \subset G_{({k_{qj}+1})}$, so $r_2^G({\underline{h}}_q)  =0$ since all the last factors in the sum defining it
are trivial. On the other hand, $r_1({\underline{h}}_q) \in G_{({l_{q1} +\ldots + l_{qm_q}})} \subset G_{({n+1})}$, so  
$R_n({\underline{h}}_1,\ldots,{\underline{h}}_p) \cdot v  =  R_n({\underline{h}}_1,\ldots,\hat{\underline{h}}_q,\ldots, {\underline{h}}_p) \cdot v$. Thus we can finally assume that for all $q$ from $1$ to $p$, 
$l_{q1} +\ldots + l_{qm_q} = n$, i.e., $l_{qj} =
k_{qj}$ for $1 \le j \le m_q$ by \REF{l1+lm}. Here  $r_2^G({\underline{h}}_q)  = [h_{q1} G_{({l_{q1}+1})}, \ldots, h_{qm_q} G_{({l_{qm_q}+1})} ] = 
[h_{q1},[ \ldots,[h_{q(m_q-1)},h_{qm_q}]\ldots] G_{({l_{q1} +\ldots + l_{qm_q}+1})} = r_1({\underline{h}}_q) G_{({n+1})}$ by definition of
the Lie ring \,L$^{\calG}(G)$. On the other hand, $r_1({\underline{h}}_1) \cdots r_1({\underline{h}}_p) G_{({n+1})} = \sum_{q=1}^p  r_1({\underline{h}}_q) G_{({n+1})}$ since 
$  r_1({\underline{h}}_q) \in G_{({l_{q1} +\ldots + l_{qm_q}})} =G_{(n)}$ for each $q$. So $R_n({\underline{h}}_1,\ldots,{\underline{h}}_p) \cdot v   =
0$, as desired.\hfbox\V


\N Proposition \ref{kerneltheta} implies that the quotient group
  \[   {\rm \bar{U}}^{\calG}(G,H) \stackrel{def}{=}  \UGH{}  \Big/ \sum_{n\ge 2} \ULG{}\,{\cal R}_n^{\calG H}  \]
is a graded $\ULG{}$--${\rm UL}(H)$-bimodule, and that
  $\theta^{\calG H}$ induces a surjective homomorphism of graded $\ULG{}$--${\rm UL}(H)$--bimodules
   \[  \sepi{\bar{\theta}^{\calG H}}{ {\rm \bar{U}}^{\calG}(G,H)\:}{\:{\rm
Gr}^{\calG \gamma}(\Z(G)I(H))}\:. \ruled\] 
Note that  $ {\rm \bar{U}}^{\gamma}(G,G) = {\rm U}^{\gamma}(G,G)\: \cong \:  {\rm \bar{U}L}(G)$ and that  
$\bar{\theta}^{\gamma G} = {\theta}^{\gamma G}$ coincides in positive degrees with the map $\theta^{\gamma}$ constructed by
Quillen. By analogy with 
the fact that $\theta^{\calG} \ot \Q$ is always an isomorphism we pose the following\vspace{3mm}

\begin{pbm}\label{conj} 
Is it true that the epimorphism 
\[ \sepi{\bar{\theta}_n^{\calG H} \ot \Q}{{\rm \bar{U}}^{\calG}_n(G,H)\ot \Q}{Q_n^{\calG}(G,H)\ot \Q} \]
 is an isomorphism for all groups
  $G$, subgroups $H$ and $n\ge 1$? In other words, is \Ker{\theta^{\calG H}} a torsion group?
\end{pbm}\vspace{2mm}

The answer is affirmative in degree $n\le 2$ by Proposition \ref{xi1et2} and also in degree 3 by Corollary \ref{conjn=3} below.


\section{The third Fox quotient}

In this section all proofs are postponed to the end.

The structure of $Q_3^{\calG}(G,H)$ is completely determined by Proposition \ref{UnGH} and the following result.\V

\begin{satz}\label{Fox3} \quad For all groups $G$ and subgroups $H$ there is a natural  exact sequence
 \[ \ruleu {\rm Tor}_1^{\ifmmode{Z\hskip -4.8pt Z} \else{\hbox{$Z\hskip -4.8pt Z$}}\fi}
 (G^{AB},   H\ab) \:\oplus\: \Big(\Ker{l_2^{\calG H}}\cap \Ker{c^H_2} \Big)
\:\mr{(\delta_1,\delta_2)}\: \UGH{3} \:\mr{\theta_3^{\calG H}} \:Q_3^{\calG}(G,H) \:\to\: 0 
\ruled \] 
Here $\delta_1$ is a homomorphism while   $\delta_2$ is a secondary operator which is a welldefined
homomorphism only modulo
\Imm{\delta_1}, i.e.\ an additive relation with indeterminacy \Imm{\delta_1}. The construction of $\delta_1$ and $\delta_2$ is given in \REF{del1def} and \REF{del2def} below in a functorial
manner and in
\REF{del1expl} and \REF{del2expl} in an explicit form.
\end{satz}\vspace{0mm}


%

%
%

\begin{kor}\label{conjn=3} Problem \ref{conj} has an affirmative answer for $n=3$.
\end{kor}

  Let us discuss Theorem \ref{Fox3} in a number of special cases. First suppose that $H=G$.
  
   If
$\calG=\gamma$  the map $l_2^{\calG G}$ is injective; hence the theorem formally generalizes
the description of
\Ker{\theta^{\calG}_3} in
\cite{Q3} (in the case $\calG =\gamma$). If $\calG \neq \gamma$ the result is still of interest as groups of
the type $Q_n^{\calG}(G,G)$ occur in the study of Fox quotients of semidirect products, see \cite{Kh12}. In particular, if $G$ is
the semidirect product  of a normal subgroup $H$ and a subgroup $K$ the quotient
  \[ X = \frac{I^3(H) \oplus I([H,K])I(H)}{I^4(H) + I(H)I([H,K])I(H) + I([H,K,H])I(H) + I([H,K,K])I(H)} \]
is proved to be a direct summand of $Q_3(G,H)$ but is not computed in \cite{Kh12}; we here fill this gap noting that  $X=Q_3^{\cal
H}(H,H)$ where the N-series 
${\cal H} = (H_{(n)})_{n\ge 1}$ is given by
$H_{(n+1)} = [H_{(n)},G]$, see also \cite{Foxsplit}.
  Indeed, the structure of $Q_3^{\calG}(G,G)$ is determined by the following two corollaries, first from a functorial point of view
further developed in remark \ref{philrem}, then by means of an explicit formula.

\begin{kor}\label{Q3(G,G)}  The group  $Q_3^{\calG}(G,G)$ is determined by the following tower of successive natural quotients
  \[  \begin{array}{rcc}
{\rm Tor}_1^{Z\!\!\!Z}(G\AB,G\ab)  &  \mr{\delta_1}  &  {\rm U}_3^{\calG}(G,G) \ruled\cr
  &  &  \surdown{q_1} \cr
{\rm Ker}\Big( \T{c}_2^{G}\,\colon\, (G_{(2)}/G_2) \sm (G_{(2)}/G_2) \to G_2/G_3 \Big)  & \mr{\delta_{21}}  &  \coker{\delta_1} \ruleu\ruled\cr
  &  &  \surdown{q_{21}} \cr
{\rm Ker}\Big( [\,,\,]\tau_1\,\colon\, {\rm Tor}_1^{Z\!\!\!Z}(G\AB,G\AB) \to G_2/[G_{(2)}\,,G_{(2)}]G_3 \Big)  & \mr{\delta_{22}} 
& 
\coker{\delta_{21}} \ruleu\ruled\cr
  &  &  \surdown{q_{22}} \cr
  &  &  \coker{\delta_{22}} \ruleu\ruled\cr
  &  &  \isodown{\tilde{\theta}_3^{\calG G}}  \cr
  &  &  Q_3^{\calG}(G,G)\ruleu\ruled
\end{array}\]
where $\T{c}_2^{G}$ is given by restriction of $c_2^G$, $\tilde{\theta}_3^{\calG G}$ is induced by ${\theta}_3^{\calG G}$, 
$\tau_1$ appears in the following part of a 6-term-exact sequence
  \[ {\rm Tor}_1^{Z\!\!\!Z}(G\AB,G\AB) \hmr{\tau_1}  G\AB \ot (G_{(2)}/G_2)  \hmr{id\ot i} G\AB \ot G\ab \hmr{id\ot \pi} G\AB\ot
G\AB
\to 0\:,\]
$[\,,\,] \,\colon\, G\AB \ot (G_{(2)}/G_2) \:\to\: G_2/[G_{(2)}\,,G_{(2)}]G_3$ is induced by the commutator pairing of $G$, and
$\delta_{21},\delta_{22}$ are induced by $\delta_2$, cf.\ the proof at the end of the section.
\end{kor}\V

It actually follows from Lemma \ref{R3=} below that $\Imm{\delta_{21}} = q_1({\cal R}_3^{\calG G})$. So if we replace ${\rm
U}_3^{\calG}(G,G)$ by $\bar{{\rm U}}_3^{\calG}(G,G)$ the above tower reduces to just the two steps involving $
\delta_{1}$ and $\delta_{22}$.

\begin{bem}\label{philrem}\rm  A similar description can be given in the general case ($H\neq G$)  by adding one additional step at
the bottom of the tower. Indeed, there is an isomorphism  $\tilde{\theta}_3^{\calG H}\,\colon\, \coker{\delta_{23}} \to
Q_3^{\calG}(G,H)$ where the construction of $\delta_{21},\delta_{22}$ resembles the one in Corollary \ref{Q3(G,G)} and where
$\delta_{23}$ looks at follows:
   \[ \begin{matrix}
{\rm Ker}\bigg(  {\rm Ker}\Big(  {\rm Tor}_1^{Z\!\!\!Z}(G/HG_{(2)}\,, H\ab) \hmr{\sigma_1} {\rm SP}^2(HG_{(2)}/G_{(2)}) \Big)
\hmr{\sigma_2}  \coker{[\,,\,]\tau_1}
\bigg)  \cr
\mapdown{\delta_{23}}  \cr
\ruleu\coker{\delta_{22}} 
\end{matrix}\]
for suitable natural maps $\sigma_1,\sigma_2,\delta_{23}$ and where ${\rm SP}^2$ denotes the second symmetric tensor power.  We
renounce to give the precise definitions and the proof as this description may be of no practical use, but we mention it
in order to illustrate our guiding philosophy: any natural construction of an abelian group associated with a nilpotent group (all
kinds of augmentation and dimension quotients, homology etc.) should be functorially expressable in terms of (generally higher order) operators
between suitable polynomial endofunctors of ${\bf Ab}$ and their derived functors, applied to appropriate abelian subquotients
of the nilpotent group in question (here {\bf Ab} denotes the category of abelian groups). 
For more examples of this structural phenomenon  see also \cite{Q3}, \cite{GoG}, \cite{D3F2}, \cite{Foxsplit}, \cite{HMP}.
\end{bem}\vspace{1mm}


\begin{kor}\label{Q3(G,G)formel} There is a natural isomorphism  $Q_3^{\calG}(G,G) \:\cong\: {\rm
U}_3^{\calG}(G,G)/(U_1+U_2)$  induced by $\theta_3^{\calG G}$   where\V

\N$\bullet$ $U_1$ is the subgroup of ${\rm
U}_3^{\calG}(G,G)$ generated by the elements \[(aG_{(2)}) \ot (b^kH_3)  -  (a^kG_{(3)}) \ot (bG_2)     + {k
\choose 2} \Big( (aG_{(2)})^2 \ot (bG_2) -   (aG_{(2)}) \ot (bG_2)^2 \Big)\]
where $a,b\in G$, $k\in \Z$ such that $a^k\in G_{(2)}$
and $b^k\in G_2$, and as usual,$(aG_{(2)})^2 =(aG_{(2)}).(aG_{(2)})$, same for $(bG_{(2)})^2$;\V

\N$\bullet$ $U_2$ is the set of elements 
\[ \sum_{q=1}^p (a_q G_{(3)}) \ot (b_qG_2) - (b_q G_{(3)}) \ot (a_qG_2)
+ \sum_{r=1}^s (c_r^{k_r}G_{(3)}) \ot (d_r G_2) -  (d_r^{k_r}G_{(3)}) \ot (c_r G_2) \]  \BE\label{delta2GGformel} \hspace{0mm} + \sum_{r=1}^s {k_r \choose
2} \Big(  (c_rG_{(2)}) . \Big((d_rG_{(2)}) - (c_rG_{(2)}) \Big) \ot (d_rG_2) \Big)
- 1 \ot (gG_4)\EE
where $a_q,b_q\in G_{(2)}$, $c_r,d_r \in G$, $k_r\in \Z$ such that $c_r^{k_r},d_r^{k_r} \in G_{(2)}$ for $1\le r\le s$ and
$g= \prod_{q=1}^p [a_q,b_q] \prod_{r=1}^s [c_r,d_r^{k_r}] \in G_3$.
 
\end{kor}
   
 The proof below shows that $U_1+U_2$   is indeed a subgroup of ${\rm U}_3^{\calG}(G,G)$.\V

\N{\bf Construction of $\delta_1$.}\quad For a group $K$ and N-series ${\cal K}$ of $K$ the isomorphisms
$\theta_n^{K,\cal K}\,\colon
{\rm U}_n{\rm L}^{\cal K}(K) 
\mr{\cong} Q_2^{\cal K}(K)$,
$n=1, 2$, provide natural exact sequences 
\BE\label{U2P2G} 0 \to {\rm U}_2\LG{} \hspace{1mm} \mr{\mu_2^{\calG}} \hspace{1mm}  P_2^{\calG}(G)  \hspace{1mm} \mr{\rho^{\calG}}
\hspace{1mm}  G^{AB} \to 0 \EE 
\BE\label{U2P2H} 0 \to {\rm U}_2{\rm L}(H)  \hspace{1mm} \mr{\mu_2^H} \hspace{1mm}  P_2(H)  \hspace{1mm}
\mr{\rho^H}
\hspace{1mm}  H^{ab} \to 0 \EE 
Tensoring   sequence \REF{U2P2G} by $H\ab = {\rm U}_1{\rm L}(H)$ and   sequence \REF{U2P2H} by
$G^{AB} = {\rm U}_1\LG{}$ gives rise to natural exact sequences
  \BE\label{deftauG}
 {\rm Tor}_1^{\Z} (G^{AB},H\ab)  \hspace{1mm} \mr{\tau_G}  \hspace{1mm} {\rm U}_2\LG{} \ot {\rm U}_1{\rm L}(H)  \hspace{1mm}
\mr{\mu_2^{\calG} \ot id} 
\hspace{1mm} P_2^{\calG}(G) \ot H^{ab}
 \hspace{1mm} \mr{\rho^{\calG} \ot id}  \hspace{1mm} G^{AB} \ot H\ab \to 0\EE
\BE\label{deftauH} 
{\rm Tor}_1^{\Z} (G^{AB},H\ab)  \hspace{1mm} \mr{\tau_H}  \hspace{1mm} {\rm U}_1\LG{} \ot {\rm U}_2{\rm L}(H)  \hspace{1mm}
\mr{id \ot \mu_2^H} 
\hspace{1mm} G^{AB} \ot  P_2(H)  
 \hspace{1mm} \mr{id \ot \rho^H}  \hspace{1mm} G^{AB} \ot H\ab \to 0\EE
\N Then define
  \BE\label{del1def} \delta_1 = \nu_{12} \tau_H - \nu_{21} \tau_G \hspace{2mm}\colon\hspace{2mm} {\rm Tor}_1^{\Z} (G^{AB},H\ab)
\hspace{2mm}\lra
\hspace{2mm} \UGH{3} \:.\EE
Note that $\delta_1$ essentially is the difference between a left and a right connecting homomorphism, kind of asymmetry phenomenon which
also induces the non-trivial torsion relations in the non-abelian tensor square and the second homology of $2$-step nilpotent groups, cf.\ \cite{GoG}.

\N To describe $\delta_1$ more explicitly let $\langle  \bar{g},k,\bar{h} \rangle $ be a typical generator of ${\rm Tor}_1^{\Z}
(G^{AB},H\ab)
$, i.e. a symbol where $k\in \Z$, $\bar{g} = gG_{(2)}\in {\rm U}_1\LG{}$, $\bar{h} = hH_2 \in {\rm U}_1{\rm L}(H)$ for $g\in G$,
$h\in H$  such that
$g^k\in G_{(2)}$ and $h^k\in H_2$, cf.\ \cite{ML} V.6. Then
   \BE\label{del1expl} \delta_1 \langle  \bar{g},k,\bar{h} \rangle   =   (\bar{g} \ot (h^kH_3)  -  (g^kG_{(3)}) \ot \bar{h}     + {k
\choose 2} ( \bar{g}^2 \ot \bar{h} -   \bar{g} \ot \bar{h}^2) \:.\EE\vspace{-1mm}
 
\N where $\bar{g}^2=\bar{g}.\bar{g}$ and 
$\bar{h}^2=\bar{h}.\bar{h}$.\V


\N{\bf Explicit formula for $\delta_2$} (the functorial construction is given in \REF{del2def} below).\quad Suppose that $H^{ab}$ is
finitely generated (the general case can be deduced from Lemmas 3.5 and 3.6 in \cite{D3F2} exactly in the same way as what
follows). Then there exists a decomposition  
$H\ab \:=\:
\bigoplus_{k=1}^r \Z/d_k\Z \cdot (h_kH_2)$ with $h_k\in H$, $d_k\in \NN{}$. Let $x\in H\ab \sm H\ab$. Then $x = \sum_{1\le i < j \le
r} a_{ij} (h_iH_2) \sm (h_jH_2)$ with $a_{ij}\in \Z$. By Lemmas 3.5 and 3.6 in \cite{D3F2} one has $x\in \Ker{l_2^{GH}}$ if and
only if for all $1\le k\le r$, $\prod_{1\le i<k} h_i^{a_{ik}}   \prod_{k <j \le n} h_j^{-a_{kj}}  \in  G_{(2)}G^{d_k}$, with
$G^d=\{g^d\,|\,g\in G\}$. Now suppose that $x\in \Ker{l_2^{GH}} \cap \Ker{c_2^H}$. Then $\gamma
\,\colon=  \prod_{1\le i < j \le r}  [h_i,h_j]^{a_{ij}} \in H_3$, and
for all $1\le k\le r$ there are $g_k\st \in
G_{(2)}$ and $g_k \in G$ such that  
$\prod_{l=1}^r h_l^{\alpha_{kl}}     = g_k\st g_k^{d_k}$ where $\alpha_{kl}= a_{lk}$ if $l<k$,  $\alpha_{kl}= 0$ if $l=k$, and
$\alpha_{kl}= - a_{kl}$ if $l>k$.
 Then

   \begin{eqnarray}\label{del2expl} \delta_2(x)   &=&  \sum_{k=1}^r \left( (g_k\st G_{(3)}) \ot  (h_kH_2)
+ (g_kG_{(2)}) \ot  (h_k^{d_k} H_3)  +  {d_k \choose 2} (g_kG_{(2)}) . \Big( (g_kG_{(2)})   \right.\nonumber\\
  &   & {} - \hspace{2mm} (h_kG_{(2)}) \Big)     \ot  (h_kH_2)
    \hspace{2mm} - \hspace{2mm}     1 \ot \bigg( (\gamma H_4) +  \sum_{l=1}^r {\alpha_{kl} \choose 2} (h_lH_2)^2 . (h_kH_2)
\nonumber\\ & & {} + \left. \phantom{{ \choose }} \sum_{1\le p < q \le r} \alpha_{kp} \alpha_{kq}  (h_pH_2) . (h_qH_2) . (h_kH_2) 
\bigg) \right) +
\Imm{\delta_1}
\end{eqnarray}

\N with $(h_lH_2)^2 = (h_lH_2).(h_lH_2)$.\V

The starting point of the proof of Theorem \ref{Fox3} is the following description of the third relative Fox polynomial group
which is completely analogous to our description of the second one in Theorem \ref{Fox2IH}. 


Let  $N$ be a   normal  subgroup of $G$ and  consider the following diagram.
 \BE\label{squFox3} \begin{matrix}
\ruleu\ruled P_2(H) \sm P_2(H)  &  \mr{l_3^{\calG H}}  &  P_2^{\bar{\cal G}}(G/N) \ot_H  P_2(H) \cr
\surdown{c_3^{H}}  &  & \surdown{\mu_3^{\calG H}} \cr
\rule{0mm}{6mm} H_2/H_4  &  \mr{d_3^{\calG H}}  &    I(G) I(H) \Big/  (\Z(G)I(N)I(H) + I_{\calG}^3(G) I(H))\rule[-5mm]{0mm}{3mm} 
\end{matrix} \rule[-5mm]{0mm}{3mm} \EE
\rule{0mm}{7mm} where for $a,b,h \in H$, $g \in G$ and $x \in H_2$
  \[ l_3^{\calG H}( p_2(a) \sm p_2(b) ) =  p_2(a) \ot p_2(b) - p_2(b) \ot p_2(a) - [p_2(a)\,,p_2(b)] \ot (p_2(a) + p_2(b)) \:,\]
 \[ c_3^{H}((aH_2) \sm (bH_2)) = [a,b]H_4\:,\hfill d_3^{\calG H}(xH_4) =  x-1 + \Z(G)I(N)I(H) + I^3_{\calG}(G) I(H)\:,\] \[
\mbox{and}\hspace{3mm} \mu_3^{\calG H}((p_2(g)\ot
(p_2(h))) =  (g-1)(h-1)  + \Z(G)I(N)I(H) + I_{\calG}^3(G) I(H)\:. \]

\N This diagram provides a simultaneous functorial description of both the third relative Fox polynomial
group and the third relative Fox subgroup, in exactly the same way as the second Fox
quotient and the second  Fox subgroup are determined in Theorem \ref{Fox2IH}, as
follows.\vspace{0.3mm}

\begin{satz}\label{Fox3push}\quad Diagram \REF{squFox3} is a pushout square of abelian
groups; in particular
  \[ \ruleu  I(G)I(H)/(\Z(G)I(N)I(H) + I_{\calG}^3(G) I(H)) \:\cong\: (P_2^{\bar{\cal G}}(G/N) \ot_H  P_2(H)) \Big/ l_3^{\calG H} ( {\rm Ker}(c_3^{H})) \]
\[  \Big(G \cap (1 + \Z(G)I(N)I(H) + I_{\calG}^3(G)I(H)) \Big) \Big/  H_4 = {\rm Ker}(d_3^{\calG H}) = c_3^{H} ({\rm Ker}(l_3^{\calG H})) \:.\]
\end{satz}\vspace{1mm}

The second equality should allow to explicitly compute the third relative Fox subgroup, in a  similar way as we deduce in \cite{D3F2} the second relative Fox-subgroup  from Theorem \ref{Fox2IH}.\V

\N{\bf Proof\,:}\quad Same principle as in the proof of \ref{Fox2IH}: diagram \REF{squFox3} is a pushout as it is obtained by gluing
the following two pushouts where ${\cal H}= \gamma_H$:
 \BE\label{squFox3proof} \begin{matrix}
\ruleu\ruled P_2(H) \sm P_2(H)  &  \mr{l_3^{{\cal H} H}}  &  P_2(H) \ot_H  P_2(H)  & \mr{i\ot id}  & P_2^{\bar{\cal G}}(G/N) \ot_H 
P_2(H)\cr
\surdown{c_3^{H}}  &   & \surdown{\mu_3^{{\cal H} H}=\mu_H} &    & \surdown{\mu_3^{\calG H}=\mu_G} \cr
\rule{0mm}{6mm} H_2/H_4  &  \mr{d_3^{{\cal H} H}}  &  I^2(H)/I^4(H)  & \mr{j}  &  \frac{\ruleu\dst I(G) I(H) }{\dst   \Z(G)I(N)I(H)
+ I_{\calG}^3(G)I(H)\ruled}
\rule[-5mm]{0mm}{3mm} \end{matrix} \rule[-5mm]{0mm}{3mm}\EE
In fact, the left hand square is a pushout by Theorem 3.5 in \cite{Q3} and the right hand square by Proposition
\ref{allgsequ}.\hfbox\V

\N{\bf Proof of Theorem \ref{Fox3}.}\hspace{2mm} First note that by \REF{U2P2G}, \REF{U2P2H} homomorphisms
  \[ \ULG{2} \ot H\ab   \hspace{2mm}\mr{\alpha_G}\hspace{2mm}  \frac{P_2^{\calG}(G) \ot P_2(H)}{\Imm{Q_2^{\calG}(G) \ot P_2(H)^2}}  
   \hspace{2mm}\ml{\alpha_H}\hspace{2mm}  G^{AB} \ot \UL{2}{H} \]
are welldefined as being factorizations of the maps 
\[ \ULG{2} \ot P_2(H)   \hspace{2mm}\mr{q(\mu_2^{\calG}\ot id)}\hspace{2mm}  \frac{P_2^{\calG}(G) \ot P_2(H)}{\Imm{Q_2^{\calG}(G) \ot P_2(H)^2}}  
   \hspace{2mm}\ml{q(id \ot \mu_2^H)}\hspace{2mm}  P_2^{\calG}(G) \ot \UL{2}{H} \]
through $id \ot \rho^H$ and $\rho^{\calG} \ot id$, resp., where $q\,\colon\, P_2^{\calG}(G) \ot P_2(H) \to \frac{P_2^{\calG}(G) \ot P_2(H)}{{\rm
Im}(Q_2^{\calG}(G) \ot P_2(H)^2)}$ is the canonical projection. Consider the following commutative diagram where $\beta_G = q_{\otimes}
\alpha_G$ with
$q_{\otimes} \,\colon\,\frac{P_2^{\calG}(G) \ot P_2(H)}{{\rm
Im}(Q_2^{\calG}(G) \ot P_2(H)^2)} \hspace{1mm}\to \hspace{1mm}  P_2^{\calG}(G) \ot_{H} P_2(H)$ being the canonical surjection.
\\
\N\makebox[14.7cm]{ \makebox[0mm]{
\begin{minipage}{20cm} 
   \[ \begin{matrix}
\ruleu\ruled P_2(H) \sm P_2(H)  &  \mr{l_3^{\calG H}}  &  P_2^{\calG}(G) \ot_H  P_2(H)  &  \ml{\beta_G}  &  \ULG{2} \ot H\ab  &=& \ULG{2}
\ot
\UL{1}{H}  \cr
\surdown{c_3^{H}}  &  & \surdown{\mu_3^{\calG H}} &  &  &  & \surdown{\nu_{21}} \cr
\rule{0mm}{6mm} H_2/H_4  &  \mr{d_3^{\calG H}}  &    \frac{\ruleu\dst I(G) I(H)}{\dst I_{\calG}^3(G)I(H)} \rule[-5mm]{0mm}{3mm} 
 & \hookleftarrow  &  Q_3^{\calG}(G,H)   &  \surltop{\theta_3^{\calG H}} & \UGH{3}
\end{matrix}\rule[-5mm]{0mm}{3mm}\]
\end{minipage}\ruled
} }\vspace{4mm}

\N As $\Ker{\mu_3^{\calG H}} = l_3^{\calG H} \Ker{c_3^{H}}$ by  Theorem \ref{Fox3push} we have 
  \BE\label{Kertheta3} \Ker{\theta_3^{\calG H}} = \nu_{21} \, \beta_G^{-1} \, l_3^{\calG H}  \, \Ker{c_3^H}\:.\EE
Thus the proof naturally divides in three steps:\V

\N{\em Step 1:}\hspace{3mm} proving the identity
  \BE\label{step1id} \nu_{21}\Ker{\beta_G}  = \Imm{\delta_1} \:; \EE
\N{\em Step 2:}\hspace{3mm}  giving a natural construction of $\delta_2$ and showing that 
  \BE\label{step2id}  \nu_{21} \, \beta_G^{-1} \, l_3^{\calG H}  \, \Ker{c_3^H}  \hspace{2mm}\equiv\hspace{2mm} \Imm{\delta_2}
\hspace{2mm}
\mbox{mod $\Imm{\delta_1}$.} \EE
\N{\em Step 3:}\hspace{3mm}  verifying the explicit formula for $\delta_2$ given by \REF{del2expl}.\V

\N{\bf Step 1.}\hspace{3mm} Let $m_2^{\cal G}\colon G^{AB} \ot G^{AB}\to {\rm U}_2\LG{} $ be the
map given by multiplication in the ring ${\rm U}_2\LG{} $, and let $m_2^H=m_2^{\gamma_H}$. 
By definition of the tensor product over $H$ there is an exact sequence 
  \[ G^{AB} \ot H\ab \ot H\ab \hspace{2mm}\mr{\psi}\hspace{2mm} \frac{P_2^{\calG}(G) \ot P_2(H)}{\Imm{Q_2^{\calG}(G) \ot P_2(H)^2}} \hmr{q_{\ot}}
P_2^{\calG}(G) \ot_H P_2(H) \hspace{1mm}\to\hspace{1mm} 0 \]
with $\psi = \alpha_G(m_2^{\cal G} \ot id) i^{\calG HH}  -  \alpha_H(id \ot m_2^H) $. Consider the following commutative diagram with exact rows
and columns by \REF{U2P2G}, \REF{U2P2H}.
  \[\begin{matrix}\ruleu\ruled
G^{AB} \ot \Ker{c_2^H}  &  \mr{id \ot l_2^H}  &  G^{AB} \ot H\ab \ot H\ab  &  \Sur{id \ot m_2^H}  &  G^{AB} \ot \UL{2}{H}\ruled  \cr
\mapdown{(m_2^{\cal G} \ot id)i^{\calG HH}(id \ot l_2^H)}  &  &  \mapdown{\psi}  &  &  \mapdown{{}-id \ot \mu_2^H}\ruled  \cr
\ruleu \ULG{2} \ot H\ab  &  \mr{\alpha_G}  &  \frac{\dst P_2^{\calG}(G) \ot P_2(H)}{\dst {\rm Im}(Q_2^{\calG}(G) \ot P_2(H)^2)}  &  \Sur{\rho^{\calG}
\ot id}  &  G^{AB}
\ot P_2(H)\ruled  \cr
\surdown{q_{GH}}  &  \mapsel{\beta_G}  &  \surdown{q_{\ot}}  &  &  \surdown{id \ot \rho^H}\ruled \cr
X  &  \mr{\bar{\alpha}_G}  &  P_2^{\calG}(G) \ot_H P_2(H)  &  \Sur{\rho^{\calG} \ot \rho^H}  &  G^{AB} \ot H\ab \ruled
\end{matrix}\]
where   
\[X = {\rm Coker}\big((m_2^{\cal G} \ot id)i^{\calG HH}(id \ot l_2^H)\big) \]
 and where $q_{GH}$ is
the canonical projection. Then by the snake lemma there is a canonical connecting homomorphism
  \BE\label{omegadef} \map{\omega}{\Ker{id \ot \mu_2^H} = \Imm{\tau_H}}{X/q_{GH}\Ker{\alpha_G}} \:,\EE
$\omega = q_Xq_{GH} \alpha^{-1}_G \psi (id \ot m_2^H)^{-1}$ with $q_X\,\colon\,X \auf X/q_{GH}\Ker{\alpha_G}$, satisfying the relation
  \BE\label{qXKeralpha}  q_X\Ker{\bar{\alpha}_G}  =  \Imm{\omega} \:.\EE
In order to compute  \Imm{\omega} we consider the two summands of $\psi$ separately. As to the first one, we have
  \BE\label{omega1} q_Xq_{GH} \alpha_G^{-1} \alpha_G(m_2^{\cal G} \ot id)i^{\calG HH}(id \ot m_2^H)^{-1} \tau_H   =  q_X \overline{i^{\calG HH}} \tau_H
\EE
where $\overline{i^{\calG HH}}$ fits into the commutative diagram
  \[\begin{matrix}\ruleu
G^{AB} \ot H\ab \ot H\ab  &  \mr{i^{\calG HH}}  &  G^{AB} \ot G^{AB} \ot H\ab  \ruled \cr
\surdown{id \ot m_2^H}  &  &  \mapdown{q_{GH}(m_2^{\cal G} \ot id)}  \ruled\cr
\ruleu G^{AB} \ot \UL{2}{H}  &  \mr{\overline{i^{\calG HH}}}  &  X\ruled
\end{matrix}\]
As to the second summand of $\psi$,
 let $\langle gG_{(2)},k,hH_2 \rangle$ be a
typical generator of Tor$^{\Z}_1(G^{AB},H\ab)$ with
$g\in G$,
$h\in H$, $k\in \Z$ such that $g^k\in G_{(2)}$ and $h^k \in H_2$. Then 
  \begin{eqnarray}\label{omega2}
\alpha_H(id \ot m_2^H)(id \ot m_2^H)^{-1} \tau_H \,\langle gG_{(2)},k,hH_2 \rangle  & = & \alpha_H\Big( (gG_{(2)}) \ot (\mu_2^H)^{-1}(kp_2(h))
\Big)  \nonumber \\
  &=&  q(id \ot \mu_2^H) \Big(p_2(g) \ot (\mu_2^H)^{-1}(kp_2(h))\Big)   \nonumber \\
  &=&  q(p_2(g) \ot k p_2(h))  \nonumber \\
  &=&  q(k p_2(g)\ot  p_2(h)) \nonumber \\
  &=&  q(\mu_2^{\calG} \ot id) \Big( (\mu_2^{\calG})^{-1}(kp_2(g)) \ot p_2(h) \Big)   \nonumber \\
  &=&  \alpha_G \Big( (\mu_2^{\calG})^{-1}(kp_2(g)) \ot (hH_2) \Big)   \nonumber \\
  &=&  \alpha_G \tau_G \,\langle gG_{(2)},k,hH_2 \rangle
\end{eqnarray}
From \REF{omega1} and \REF{omega2} we  obtain the relations
  \[  \omega \tau_H  =  q_X \left( \overline{i^{\calG HH}} \tau_H - q_{GH} \tau_G \right)\:,\]
whence by \REF{qXKeralpha}
  \BE\label{Keralpha}  \Ker{\bar{\alpha}_G}  =  \Imm{\overline{i^{\calG HH}} \tau_H - q_{GH} \tau_G}  +  q_{GH} \Ker{\alpha_G}\:.\EE
Next we compute \Ker{\alpha_G}. Consider the following commutative diagram whose top row is part of the six-term exact sequence
obtained by tensoring sequence \REF{U2P2G} with $P_2(H)$.
  \[\begin{matrix}\ruleu
{\rm Tor}_1^{\ifmmode{Z\hskip -4.8pt Z} \else{\hbox{$Z\hskip -4.8pt Z$}}\fi}(G^{AB},P_2(H))  &  \mr{\tau}  & \ULG{2} \ot P_2(H)  & 
\mr{\mu_2^{\calG} \ot id}  &  P_2^{\calG}(G) \ot P_2(H)  \ruled\cr
\mapdown{{\rm Tor}_1^{\ifmmode{Z\hskip -4.8pt Z} \else{\hbox{$Z\hskip -4.8pt Z$}}\fi}(id,\rho^H)}  &  &  \surdown{id \ot \rho^H}  & 
&  \surdown{q}    \ruled\cr
\ruleu {\rm Tor}_1^{\ifmmode{Z\hskip -4.8pt Z} \else{\hbox{$Z\hskip -4.8pt Z$}}\fi}(G^{AB},H\ab)  &    \mr{\tau_G}  & \ULG{2} \ot
H\ab  &  \mr{\alpha_G}  &   \frac{\dst P_2^{\calG}(G) \ot P_2(H)}{\dst {\rm Im}(Q_2^{\calG}(G) \ot P_2(H)^2)} \ruled
\end{matrix}\]
The left hand square commutes by naturality of the connecting homomorphism and the right hand square by definition of $\alpha_G$.
As $\Ker{q} = (\mu_2^{\calG}\ot id)\Ker{id \ot \rho^H}$ we have
   \begin{eqnarray}
q_{GH}\Ker{\alpha_G}  &=&  \Imm{q_{GH}\, \tau_G {\rm Tor}_1^{\ifmmode{Z\hskip -4.8pt Z} \else{\hbox{$Z\hskip -4.8pt
Z$}}\fi}(id,\rho^H) } \nonumber\\
\label{Kerbeta}  &=&   {\rm Im}\left( (\overline{i^{\calG HH}} \tau_H - q_{GH} \tau_G ) 
\rond {\rm Tor}_1^{\ifmmode{Z\hskip -4.8pt Z} \else{\hbox{$Z\hskip -4.8pt Z$}}\fi}(id,\rho^H)   \right)
\end{eqnarray}
since $\tau_H \rond {\rm Tor}_1^{\ifmmode{Z\hskip -4.8pt Z} \else{\hbox{$Z\hskip -4.8pt Z$}}\fi}(id,\rho^H) = 0$ as these are
consecutive maps in the six term exact sequence part of which is displayed in \REF{deftauH}. Thus by \REF{Keralpha},
  \BE\label{Keralphafinal} \Ker{\bar{\alpha}_G}  =  \Imm{ \overline{i^{\calG HH}} \tau_H - q_{GH} \tau_G } \:.\EE
Now we are ready to compute $\nu_{21}\Ker{\beta_G}$. First of all we note that 
  \[ \nu_{21}(m_2^{\cal G} \ot id)i^{\calG HH} (G^{AB} \ot l_2^H \Ker{c_2^H} )  =  0 \]
since 
  \BE\label{nu21=nu12}  \nu_{21}(m_2^{\cal G} \ot id)i^{\calG HH} = \nu_{12}(id \ot m_2^H)  \EE
by middle associativity of the tensor product over $\UL{}{H}$ used to define $\UGH{}$. Hence $\nu_{21}$ factors as $\ULG{2}
\ot H\ab  \Sur{q_{GH}} X \mr{\bar{\nu}_{21}} \UGH{3}$, and $\nu_{21} \Ker{\beta_G} = \bar{\nu}_{21} \Ker{\bar{\alpha}_G}$. But
$\bar{\nu}_{21} 
\overline{i^{\calG HH}}  =
\nu_{12}$ since 
   \begin{eqnarray*}
\bar{\nu}_{21}  \overline{i^{\calG HH}} (id \ot m_2^H)  &=&  \bar{\nu}_{21} q_{GH} (m_2^{\cal G} \ot id) i^{\calG HH} \\
  &=&   {\nu}_{21}  (m_2^{\cal G} \ot id) i^{\calG HH} \\
  &=&   {\nu}_{12}   (id \ot m_2^H) 
\end{eqnarray*}
again by \REF{nu21=nu12}. Hence 
  \[\nu_{21} \Ker{\beta_G}  =  \Imm{\nu_{12}\tau_H  -  \nu_{21}\tau_G }  = \Imm{\delta_1}  \:,\]
as desired.\V

\N{\bf Step 2.}\hspace{3mm}  Consider the following diagram.
\vspace{2mm}

\N\makebox[14.7cm]{ \makebox[0mm]{
\begin{minipage}{18cm}\small
  \BE\label{del2dia}\begin{matrix}
\ruleu\ruled  \ULG{2} \ot \UL{1}{H}  &  =  &  \ULG{2} \ot H\ab  &  \mr{\beta_G}  &  P_2^{\calG}(G) \ot_H P_2(H)  &    
\Sur{\rho^{\calG} \ot
\rho^H}  &  G^{AB} \ot H\ab  \cr
\mapdown{\nu_{21}}  &  &  \mapup{l_3^{\,\prime}}  &  & \mapup{l_3^{\calG H}}  &  &  \mapup{l_2^{\calG H}} \ruled\cr
\ruleu \UGH{3}  &  &  H\ab \ot H\ab \ot P_2(H)  &  \mr{\lambda}  &  P_2(H) \sm  P_2(H)  &  \Sur{\rho^H \sm
\rho^H}  &  H\ab \sm H\ab  \ruled\cr
\mapup{\nu_{03}}  &  &  \surdown{c_3^{\,\prime}}  &  & \surdown{c_3^{H}}  &  & \surdown{c_2^{H}}  \ruled\cr
\ruleu 1_{\ULG{}} \ot \UL{3}{H}  &  \ml{i_3}  &  H_3/H_4  &  \stackrel{j_3}{\hra}  &  H_2/H_4  &  \Sur{}  &  H_2/H_3 \ruled
\end{matrix}\EE
\end{minipage}\ruled
} }\vspace{6mm}

\N Here $l_3\st = (m_2^{\cal G} \ot id) i^{HHH} l_{33}^H (id \ot id \ot \rho^H)$, $c_3\st = c_{33}^H (id \ot id \ot \rho^H)$, $\lambda =
q_{\sst\sm} (\mu_2^H m_2^H \ot id)$ with $q_{\sst\sm}\,\colon\,P_2(H) \ot  P_2(H)  \auf  P_2(H) \sm P_2(H)$, and
$i_3\,\colon\, H_3/H_4 ={\rm L}_3(H)  \to \UL{3}{H}  \to 1_{\ULG{}} \ot \UL{3}{H}$ is the composite of the canonical maps.

Diagram \REF{del2dia} commutes; this is clear from the definitions for the two rightmost squares and was essentially proved in
\cite[Lemma 4.3]{Q3}  for the two middle squares. For  the left hand rectangle this is due to middle associativity of the tensor
product over $\UL{}{H}$ used to define $\UGH{}$. Moreover, omitting the left hand rectangle, the rows of the diagram are exact by
\REF{U2P2G}, \REF{U2P2H}. Now define the additive relation $\delta_2\,\colon\, \Ker{c_2^{H}} \,\cap\,\Ker{l_2^{\calG H}}  \to \UGH{3}$ by
  \BE\label{del2def}  \delta_2  =  (\nu_{21} \beta_G^{-1} l_3^{\calG H}  -  \nu_{03} i_3j_3^{-1}c_3^H)(\rho^H \sm \rho^H)^{-1} \:.\EE
The indeterminacy of the first factor from the left is $\nu_{21}\Ker{\beta_G} = \Imm{\delta_1}$ by \REF{step1id}, and the 
indeterminacy of  the second factor is annihilated by the first modulo \Imm{\delta_1} by commutativity of diagram \REF{del2dia} and
exactness of its middle row. Whence $\delta_2$ is a welldefined homomorphism modulo \Imm{\delta_1}.

Let $\sepi{q_3}{\UGH{3}}{\UGH{3}/\Imm{\delta_1}}$ and define the homomorphism 
  \[\delta_2\st  \,\colon\hspace{2mm}  \Ker{c_3^H}  \hspace{1mm}\cap\hspace{1mm} (l_3^{\calG H})^{-1}\Imm{\beta_G} \hspace{2mm}
\lra\hspace{2mm} \UGH{3}/\Imm{\delta_1}\:,\hspace{3mm}  \delta_2\st = q_3 \nu_{21} \beta_G^{-1} l_3^{\calG H} \:.\]
By \REF{Kertheta3} one has $q_3 \Ker{\theta_3^{\calG H}}  =  \Imm{\delta_2\st}$. Now by the snake lemma, diagram \REF{del2dia} induces an
exact sequence 
  \[ \Ker{c_3\st}  \mr{} \Ker{c_3^H}  \mr{} \Ker{c_2^H}   \mr{} \coker{c_3\st} =  0 \]
which implies another exact sequence
  \BE\label{Kerc3sequ}  \Ker{c_3\st}  \mr{\lambda\st}  \Ker{c_3^H}  \hspace{1mm}\cap\hspace{1mm} (l_3^{\calG H})^{-1}\Imm{\beta_G}
\hspace{2mm} \mr{\rho^{\,\prime}} \hspace{2mm} \Ker{c_2^{H}} \,\cap\,\Ker{l_2^{\calG H}} \:\to\:0\EE
where $\lambda\st$ and $\rho\st$ are the restrictions of  $\lambda$ and $\rho^H \sm \rho^H$, resp.
Now $\delta_2\st \lambda\st  =  q_3 \nu_{21} l_3\st  =  q_3 \nu_{03} i_3 c_3\st = 0$, so $\delta_2\st$ factors as 
  \[
\Ker{c_3^H}  \hspace{1mm}\cap\hspace{1mm} (l_3^{\calG H})^{-1}\Imm{\beta_G} \hspace{2mm}  \Sur{\rho\st} 
 \Ker{c_2^{H}} \,\cap\,\Ker{l_2^{\calG H}}  \hmr{\bar{\delta}_2\st}     \UGH{3}/\Imm{\delta_1}\:.\] Let $x\in \Ker{c_2^{H}}
\,\cap\,\Ker{l_2^{\calG H}}$.
 By \REF{Kerc3sequ} there is $x\st \in \Ker{c_3^H}  \hspace{1mm}\cap\hspace{1mm} (l_3^{\calG H})^{-1}\Imm{\beta_G}$ such that
$\rho\st(x\st) = x$. Then $\bar{\delta}_2\st(x) = \delta_2\st(x\st)  =  q_3(\nu_{21} \beta_G^{-1}l_3^{\calG H}  -  \nu_{03} i_3
j_3^{-1}c_3^H)(x\st)  =  q_3\delta_2(x)$. Thus $\bar{\delta}_2\st  = q_3\delta_2$ and $q_3\Ker{\theta_3^{\calG H}} = \Imm{\delta_2\st} 
= 
\Imm{\bar{\delta}_2\st} = q_3\Imm{\delta_2}$ which proves relation \REF{step2id}.\V

\N{\bf Step 3.}\hspace{3mm}  The explicit formula for $\delta_2$ in \REF{del2expl} is obtained by a straightforward calculation,
taking $\sum_{1\le i < j \le r} a_{ij} p_2(h_i) \sm p_2(h_j) $ as a representative element of the coset $(\rho^H \sm
\rho^H)^{-1}(x)$, and using the following device: for $a,b\in G$, $p_2(ab) = p_2(a) + p_2(b)  +  p_2(a)  p_2(b)$; in particular, if
$a$ or $b$ is in $G_{(2)}$, $p_2(ab) = p_2(a) + p_2(b)$. More generally, for $x_1,\ldots,x_n \in G$, $p_2(\prod_{i=1}^n x_i)  = 
\sum_{i=1}^n p_2(x_i)  +  \sum_{1\le i < j \le n} p_2(x_i) p_2(x_j)$. In particular, $p_2(x^n) = np_2(x) + {n \choose 2} p_2(x)^2$
which is also true for negative $n$. This implies the equation 
  \[ d_k p_2(g_k) \ot p_2(h_k)  = p_2(g_k) \ot \left( p_2(h_k^{d_k}) - {d_k \choose 2} p_k(h_k)^2 \right)\:.\]
 Finally, for $g\in G$
and $h\st \in H_2$,
  \[  \nu_{21} \beta_G^{-1} (p_2(g ) \ot p_2(h\st))  =  \nu_{12}((gG_{(2)}) \ot (h\st H_3))  +  \Imm{\delta_1}  \]
which is deduced from the following relations for $h_1,h_2 \in H$, $h_i\st\in H_2$:
\begin{eqnarray*}
p_2({\textstyle\prod_i} h_i\st) &=& {\textstyle\sum_i}p_2(h_i\st)\\
p_2([h_1,h_2]) &=& [p_2(h_1),p_2(h_2)]\\
\nu_{21} \beta_G^{-1}(p_2(g)\ot p_2(h_1)p_2(h_2)) 
&= & \nu_{21} \beta_G^{-1}(p_2(g)p_2(h_1)\ot p_2(h_2)) \\
&= &(gG_{(2)}) . (h_1G_{(2)}) \ot (h_2H_2)   +  \Imm{\delta_1}  \\
  &=&  (gG_{(2)}) \ot (h_1 H_2) . (h_2H_2)   +  \Imm{\delta_1}  \\
  &=&  \nu_{12}((gG_{(2)}) \ot (h_1H_2) . (h_2H_2) ) +  \Imm{\delta_1} 
  \end{eqnarray*} 
\hfbox\V

\N We now establish the link between the description of \Ker{\theta_3^{\calG H}} given by Theorem \ref{Fox3} and the subgroup ${\cal R}_3^{\calG H}$ of \Ker{\theta_3^{\calG H}} constructed in section 5. First note that \Ker{l_2^{\calG H}} contains the canonical subgroup \[ \Gamma
={\rm Im}\Big( \big((H\cap G_{(2)})/H_2\big) \sm \big( (H\cap G_{(2)})/H_2\big) \hspace{1mm}\to{}\hspace{1mm} (H/H_2) \sm (H/H_2)\Big) \:.\]

\begin{lem}\label{R3=} One has   $\delta_2(\Gamma \cap \Ker{c_2}) \equiv {\cal R}_3^{\calG H}$ mod \Imm{\delta_1}.
\end{lem}

\proof  Let $x\in \Gamma \cap \Ker{c_2}$. Then $x= \sum_{q=1}^p (h_qH_2) \sm (h_q\st H_2)$ with $h_q,h_q\st \in H \cap G_{(2)}$,
$q=1,\ldots,p$, such that $h \colon= \prod_{q=1}^p [h_q,h_q\st] \in H_3$. First of all, note that $x\st = 
 \sum_{q=1}^p   p_2(h_q) \sm p_2(h_q\st) \in (\rho^H \sm \rho^H)^{-1} (x)$. Then using the fact that
 $[ p_2(h_q)\,,p_2(h_q\st)] = 0$ in $P_2^{\calG}(G)$
since $ p_2(h_q) \in p_2(G_{(2)}) \subset Q_2^{\calG}(G)$ we get
  \begin{eqnarray*}
l_3^{\calG H}(x\st)  &=& \sum_{q=1}^p \Big( p_2(h_q) \ot p_2(h_q\st) - p_2(h_q\st) \ot p_2(h_q) \Big) \\
  &=& \beta_G \bigg( \sum_{q=1}^p (h_qG_{(3)}) \ot  (h_q\st H_2)  - (h_q\st G_{(3)}) \ot  (h_q H_2) \bigg)  \:.
  \end{eqnarray*}
 On the other hand,
  \[  \nu_{03}i_3j_3^{-1}c_3^H(\rho^H \sm \rho^H)^{-1} (x) \equiv \nu_{03}i_3j_3^{-1} \bigg( \prod_{q=1}^p [h_q,h_q\st] H_4 \bigg) =
\nu_{03}i_3(hH_4) = 1 \ot (hH_4) \:.\]
Therefore
  \begin{eqnarray*} 
\delta_2(x)  &=&  \sum_{q=1}^p \bigg( (h_qG_{(3)}) \ot  (h_q\st H_2)  - (h_q\st G_{(3)}) \ot 
(h_q H_2) \bigg)  - 1 \ot (hH_4)  + \Imm{\delta_1}  \\
  &=& {} - R_3({\underline{h}}_1,\ldots,{\underline{h}}_p) + \Imm{\delta_1} 
  \end{eqnarray*}
where each ${\underline{h}}_q = (h_q,h_q\st)$ is of height $\ge n=3$, with $k =1$ and $l =2$ for $h_q$ and for $h_q\st$. The rest should now be clear.
\hfbox\V

\N{\bf Proof of Corollary \ref{conjn=3}:}\hspace{2mm}   By Theorem \ref{Fox3}, $\Ker{\bar{\theta}_3^{{\cal G},H}} = q
\Imm{\delta_1,\delta_2}$ where
$q\,\colon\,\UGH{3}
\auf
{\rm \bar{U}}_3(G,H)$ is the canonical projection.
By Lemma \ref{R3=}, $\delta_2$ induces a homomorphism $\bar{\delta}_2\,\colon\, \frac{{\rm Ker}(l_2^{\calG H})\cap {\rm
Ker}(c^H_2)}{\Gamma
\cap {\rm Ker}(c^H_2)} \hspace{1mm}\lra\hspace{1mm} \frac{{\rm \bar{U}}_3(G,H)}{q {\rm Im}(\delta_1)}$ and we have an exact sequence 
$0\to q
\Imm{\delta_1}
\:\to\: q
\Imm{\delta_1,\delta_2} \:\to\:\Imm{\bar{\delta}_2}  \to 0$. But
 the quotient $\Ker{l_2^{\calG H}}/\Gamma$ is torsion  by \cite[Lemma 2.7]{D3F2}. Hence
its subgroup
$\frac{{\rm Ker}(l_2^{\calG H})\cap {\rm Ker}(c^H_2)}{\Gamma \cap {\rm Ker}(c^H_2)}$ is torsion, too. Thus being an extension of
torsion groups,
$\Ker{\bar{\theta}_3^{\calG H}}$ is torsion  which implies the assertion.
\hfbox\vspace{4mm}

\N{\bf Proof of corollaries \ref{Q3(G,G)} and \ref{Q3(G,G)formel}\,:} \quad One has the following sequence of homomorphisms
\\
\N\makebox[14.7cm]{ \makebox[0mm]{
\begin{minipage}{20cm}\small
  \[ {\rm Tor}_1^{Z\!\!\!Z}(G\AB,G\AB)  \hmr{\tau_1}  G\AB \ot (G_{(2)}/G_2) \hmr{\bar{\nu}}  \frac{G\ab \sm
G\ab}{\Imm{(G_{(2)}/G_2)
\sm  (G_{(2)}/G_2)}}  \hmr{\overline{l_2^{\calG G}}}  G\AB \ot G\ab \]
\end{minipage}\ruled
} }\vspace{4mm}

\N where $\bar{\nu}(gG_{(2)} \ot g\st G_2) = \overline{gG_2 \sm
g\st G_2}$, $g,g\st\in G$.  Now
consider the following commutative diagram where $\hat{c}_2^G,\check{c}_2^G$ are given by restriction of ${c}_2^G$.
  \[\begin{matrix}
\ruleu\Ker{\T{c}_2^G}  &  \mr{\delta_{21}}  &  {\rm U}_3^{\calG}(G,G)/\Imm{\delta_1}  &  \Sur{q_{21}}  & \coker{\delta_{21}} \cr
\|  &  &  \mapup{q_1\delta_2}  &  &    \mapup{\overline{q_1\delta_2}} \cr
\ruled\Ker{\hat{c}_2^G}  &  \hra  &  \Ker{c_2^G} \cap \Ker{l_2^{\calG G}}  &  \mr{}  &  \Ker{\overline{\check{c}_2^G}}
\cr
 \mapdown{}  &  &  \mapdown{} &  & \mapdown{}  \cr
\ruled (G_{(2)}/G_2) \sm  (G_{(2)}/G_2) &  \mr{i}  &  \Ker{l_2^{\calG G}}  &  \Sur{}  &  \Ker{\overline{l_2^{\calG G}}} \cr
\surdown{\hat{c}_2^G}  &  &  \mapdown{\check{c}_2^G}  &  &  \mapdown{\overline{\check{c}_2^G}} \cr
\ruleu\rule[-3mm]{0mm}{0mm} [G_{(2)}\,,G_{(2)}]G_3 \Big/ G_3  & \hra  &  G_2/G_3  &  \Sur{}  & G_2 \Big/ [G_{(2)}\,,G_{(2)}]G_3
\end{matrix}\]
The columns and the two bottom rows are exact, so by the snake lemma the second row is short exact as $\hat{c}_2^G$ is surjective. 
Now by Lemma 2.7 in \cite{D3F2} one has $\Ker{\overline{l_2^{\calG G}}}  = \Imm{\bar{\nu}\tau_1}$, so
$\Ker{\overline{\check{c}_2^G}} = \bar{\nu}\tau_1 \Ker{\overline{\check{c}_2^G}\bar{\nu}\tau_1} = \bar{\nu}\tau_1
\Ker{[\,,\,]\tau_1}$. So letting $\delta_{22}$ be the restriction of $\overline{q_1\delta_2} \bar{\nu}\tau_1 $ to
${\rm Ker}([\,,\,]\tau_1)$ we have $\Imm{\delta_{22}} = \Imm{\overline{q_1\delta_2}} =
q_{21} \Imm{q_1\delta_2}  =  q_{21} q_1 \Ker{\theta_3^{\calG G}}$ by
Theorem
\ref{Fox3}, whence Corollary \ref{Q3(G,G)} is proved. As to Corollary \ref{Q3(G,G)formel} first note that by \REF{del1expl},
$U_1=\Imm{\delta_1}$. Now let $x\in \Ker{l_2^{\calG G}}$. By the identity $\Ker{\overline{l_2^{\calG G}}}  = \Imm{\bar{\nu}\tau_1}$
above there is  $y\in {\rm Tor}_1^{Z\!\!\!Z}(G\AB,G\AB)$ such that $\bar{x}=\bar{\nu}\tau_1(y)$ in $G\ab
\sm G\ab/\Imm{(G_{(2)}/G_2) \sm  (G_{(2)}/G_2)}$. By \cite{ML} V.6, $y=
\sum_{r=1}^s \langle c_rG_{(2)}, k_r, d_rG_{(2)} \rangle$ with $s\ge 1$, $c_r,d_r \in G$, $k_r\in \Z$ such that
$c_r^{k_r},d_r^{k_r} \in G_{(2)}$, and $\bar{\nu}\tau_1(y) = \sum_{r=1}^s \overline{(c_r G_2) \sm (d_r^{k_r} G_2)}$. Thus $x=
\sum_{q=1}^p (a_q G_2) \sm (b_qG_2) + \sum_{r=1}^s  (c_r G_2) \sm (d_r^{k_r} G_2)$ with $q\ge 1$ and $a_q,b_q\in G_{(2)}$. Now
suppose that $x\in \Ker{c_2^G}$ which means that  $g= \prod_{q=1}^p [a_q,b_q] \prod_{r=1}^s [c_r,d_r^{k_r}] \in G_3$. To
compute $\delta_2(x)$ first note that putting $x_1 = \sum_{q=1}^p p_2(a_q) \sm p_2(b_q)$ and $x_2=\sum_{r=1}^s  p_2(c_r) \sm
p_2(d_r^{k_r})$ one has $x_1+x_2\in
(\rho^G\sm \rho^G)^{-1}(x)$. By the calculation in the proof of Lemma \ref{R3=},
\BE\label{delta21}  l_3^{\calG G} (x_1)
 =  \beta_G \bigg( \sum_{q=1}^p (a_q G_{(3)}) \ot (b_qG_2) - (b_q G_{(3)}) \ot (a_qG_2) \bigg) \:. \EE
Now note that $[p_2(c_r),p_2(d_r^{k_r})] = p_2([c_r,d_r^{k_r}])=0$ in $P_2^{\calG}(G)$ since $
[c_r,d_r^{k_r}]\in G_{(3)}$.
This  justifies the first of the following identities.

  \begin{eqnarray}
l_3^{\calG G} (x_2)
&=&  \sum_{r=1}^s  \Big(p_2(c_r) \ot p_2(d_r^{k_r})  -  p_2(d_r^{k_r}) \ot p_2(c_r) \Big) \nonumber\\
&=&  \sum_{r=1}^s  \Big( p_2(c_r) \ot \Big(k_r p_2(d_r) +   {k_r \choose 2} p_2(d_r)^2 \Big) - p_2(d_r^{k_r}) \ot p_2(c_r) \Big)
\nonumber\\ 
&=&  \sum_{r=1}^s  \Big( k_rp_2(c_r) \ot   p_2(d_r) +   {k_r \choose 2}  p_2(c_r) \ot p_2(d_r)^2 - p_2(d_r^{k_r}) \ot
p_2(c_r) \Big)
\nonumber\\
&=&  \sum_{r=1}^s  \Big( \Big(p_2(c_r^{k_r}) -    {k_r \choose 2} p_2(c_r)^2\Big) \ot p_2(d_r) +   {k_r \choose 2} 
p_2(c_r)p_2(d_r) \ot p_2(d_r) \nonumber\\ 
& &  - p_2(d_r^{k_r}) \ot p_2(c_r) \Big)
\nonumber\\
&=&  \sum_{r=1}^s  \Big( p_2(c_r^{k_r}) \ot p_2(d_r) - p_2(d_r^{k_r}) \ot p_2(c_r)  \nonumber\\ & & +   {k_r \choose 2}
p_2(c_r)
\Big(p_2(d_r)-p_2(c_r) \Big)
 \ot p_2(d_r) \Big) \nonumber\\
&=&\label{delta22} \beta_G\bigg( \sum_{r=1}^s  \Big( (c_r^{k_r}G_{(3)}) \ot (d_r G_2) -  (d_r^{k_r}G_{(3)}) \ot (c_r G_2)
  \nonumber\\ & &  +   {k_r \choose 2} \Big(  (c_rG_{(2)}) \Big((d_rG_{(2)}) - (c_rG_{(2)}) \Big) \ot (d_rG_2)  \Big)\Big)
\bigg)
\end{eqnarray}
It now follows from \REF{delta21} and \REF{delta22} that $\delta_2(x)$ is represented modulo \Imm{\delta_1} by the element given in
\REF{delta2GGformel} which achieves the proof.\hfbox\V

\begin{center} \sc Acknowledgments  \vspace{3mm} \end{center} 

{\small The author is indebted to P. Littelmann for a helpful discussion about higher commutators. Furthermore, 
he wishes to express his gratitude to the Institut de Recherche Math\'ematique Avanc\'ee, Strasbourg, and in particular to his
former director, Jean-Louis Loday, for the warm hospitality and ideal working conditions provided during the author's postdoctoral
fellowship  in 1993-1995 where part of the  research presented here was done.}

\vspace{8mm}

\tw